\documentclass[12pt,a4paper]{article}
\usepackage{amssymb,amsmath,latexsym,theorem,natbib, easybmat,
  verbatim,graphicx,color,epsfig}

\newcommand{\tr}{\mathrm{tr}}
\newcommand{\mse}{\mathrm{MSE}}
\newcommand{\eff}{\mathrm{Eff}}

\numberwithin{equation}{section}

\newcommand{\proofend}{\hfill$\square$}

\newtheorem{thm}{Theorem}[section]

\theorembodyfont{\rm}

\newtheorem{ex}[thm]{Example}

\begin{document}
\title{K-optimal designs for parameters of shifted Ornstein-Uhlenbeck
  processes and sheets}

\author{S\'andor Baran \\ 
Faculty of Informatics,
University of Debrecen, Hungary}

\date{}
\maketitle

\begin{abstract}
Continuous random processes and fields are regularly applied to model temporal or spatial phenomena in many different fields of science, and model fitting is usually done with the help of data obtained by observing the given process at various time points or spatial locations. In these practical applications sampling designs which are optimal in some sense are of great importance. 
We investigate the properties of the recently introduced K-optimal design for temporal and spatial linear regression models driven by Ornstein-Uhlenbeck processes and sheets, respectively, and highlight the differences compared with the classical D-optimal sampling. A simulation study displays the superiority of the K-optimal design for large parameter values of the driving random process.

\bigskip

\end{abstract}

\section{Introduction}
   \label{sec:sec1}

Continuous random processes and fields are regularly applied to model temporal or spatial phenomena in many different fields of science such as agriculture, chemistry, econometrics, finance, geology or physics. Model fitting is usually done with the help of data obtained by observing the given process at various time points or spatial locations. These observations are either used for parameter estimation or for prediction. However, the results highly depend on the choice of the data collection points. Starting with the fundamental works of \citet{hoel58} and \citet{kiefer59}, a lot of work has been done in the field of optimal design. Here by a design we mean a set \ $\boldsymbol\xi = \{x_1, x_2 , \ldots , x_n \}$ \ of distinct time points or locations where the investigated process is observed, whereas optimality refers to some prespecified criterion \citep{muller}. In case of prediction, one can use, e.g., the Integrated Mean Square Prediction Error criterion, which minimizes a functional of the error of the kriging predictor \citep{baz,bss13} or maximize the entropy of observations \citep{swin}. In parameter estimation problems, 
a popular approach is to consider information based criteria. An A-optimal design minimizes the trace of the inverse of the Fisher information matrix (FIM) on the unknown parameters, whereas E-, T- and D- optimal designs maximize the smallest eigenvalue, the trace and the determinant of the FIM, respectively 
\citep[see, e.g.,][]{puk,aw,pazman}. The latter design criterion for regression experiments has been studied by several authors both in uncorrelated \citep[see, e.g.,][]{silvey} and in correlated setups  \citep{ms,ks,zba09,dpz}. However, there are several situations when D-optimal designs do not exist, for instance, if one has to estimate the covariance parameter(s) of an Ornstein-Uhlenbeck (OU) process \citep{zba09} or sheet \citep{bsschemo}. This deficiency can obviously be corrected by choosing a more appropriate design criterion. In case of regression models 
a recently introduced approach, which optimizes the condition number of the FIM, called K-optimal design \citep{yz}, might be a reasonable choice. K-optimal designs try to minimize the error sensitivity of experimental measurements \citep{myz} resulting in more reliable least squares estimates of the parameters. However, one can also consider the condition number of the FIM as a measure of collinearity \citep{rz}, thus minimizing the condition number avoids multicollinearity. 

In contrast to the standard information based design criteria, the condition number (and the corresponding optimization problem) is not convex, only quasiconvexity holds \citep{myz}.  Hence, finding a K-optimal design usually requires non-smooth algorithms. \citet{yz} consider polynomial regression models and solve the K-optimal design problem with nonlinear programming, whereas in \citet{rz} simulated annealing is applied. In this class of models K-optimal designs are quite similar to their A-optimal counterparts. Further, \citet{myz} investigate Chebyshev polynomial models and suggest a two-step approach to find a probability distribution approximating the K-optimal design. 

Further, one should also mention that K-optimal design is invariant to the multiplication of the FIM by a scalar, so it does not measure the amount of information on the unknown parameters. Besides this, K-optimality obviously does not have meaning for one-parameter models, but in this case multicollinearity does not appear either.

All regression models where K-optimality has been investigated so far consider  uncorrelated errors, but there are no results for correlated processes.
In the present paper we derive K-optimal designs for estimating the regression parameters of simple temporal and spatial linear models driven by OU processes and sheets, respectively, and compare the obtained sampling schemes with the corresponding D-optimal designs. Both increasing domain and infill equidistant designs are investigated and the key differences between the two approaches are highlighted. Our aim is to give a first insight into the behaviour of K-optimal designs in a correlated setup, but many results presented here can be generalized to models with different base functions and/or correlation structures \citep[see, e.g.,][]{nather,dpz2}. This is a natural direction for further research.

\section{Ornstein-Uhlenbeck processes with linear trend}
   \label{sec:sec2}
Consider the stochastic process
\begin{equation}
   \label{model}
Y(s) = \alpha_0 + \alpha_1 s +U(s)
\end{equation}
with  design points taken from a compact  interval \ $[a,b]\subset
{\mathbb R}$, \ where $U(s), \ s\in {\mathbb R}$, is a stationary OU
process, that is a zero mean Gaussian process with covariance structure
\begin{equation}
   \label{oucov}
{\mathsf E}\,U(s)U(t)=\frac{{\sigma}^2}{2\beta}\exp\big
(-\beta|s-t|\big ),
\end{equation}
with \ $\beta>0,  \ \sigma>0$. \ We remark that
\ $U(s)$ \ can also be represented as
\begin{equation}
  \label{ourep}
U(s)=\frac{\sigma}{\sqrt{2\beta}}{\mathrm
  e}^{-\beta s}{\mathcal W}\big({\mathrm e}^{2\beta s}\big),
\end{equation}
where \ ${\mathcal W}(s), \ s\in {\mathbb R}$, \ is a  standard
Brownian motion \citep[see, e.g.,][]{sw,bpz}. In the present study
the parameters \ $\beta$ \ and $\sigma$ \ of the driving OU
process \ $U$ \ are assumed to be known. However, a valuable
direction for future research will be the investigation of models
where these parameters should also be estimated. We remark that the same type  of regression model appears in \citet{ms}, where the properties of  D-optimal design under a different driving process are investigated. 

For model \eqref{model}, the FIM \ 
$\mathcal I_{\alpha_0,\alpha_1}(n)$ \ on the unknown parameters \ $\alpha_0$ \
and \ $\alpha_1$ \ based on observations \ $\big\{ Y(s_i), \
i=1,2,\ldots ,n\big\}, \ n\geq 2,$ \ equals  
\begin{equation*}
\mathcal I_{\alpha_0,\alpha_1}(n)=H(n)C(n)^{-1}H(n)^{\top}, \quad \text{where}
\quad H(n):=\begin{bmatrix} 1 & 1 & \cdots & 1 \\ s_1 & s_2 & \cdots &
  s_n \end{bmatrix},
\end{equation*}
and  \ $C(n)$ \ is the covariance matrix of the observations
\citep[see, e.g.,][]{xia,pazman}. Without loss of generality, one can set the variance of \ $U$ \ to be equal to one, which reduces  \ $C(n)$ \ to a correlation matrix. Due to the particular structure of \ $C(n)$
\ resulting in a special form of its inverse (see \ref{subs:subsA.1} or \citet{ks}), a short
calculation shows that 
\begin{equation*}
\mathcal I_{\alpha_0,\alpha_1}(n)=
\begin{bmatrix} L_1(n) & L_2(n) \\ L_2(n) & L_3(n) \end{bmatrix},
\end{equation*}
with
\begin{equation}
  \label{entry}
L_1(n):=\!1\!+\!\sum_{i=1}^{n-1} \frac{1\!-\!p_i}{1\!+\!p_i}, \qquad \
L_2(n):=\!s_1\!+\!\sum_{i=1}^{n-1} 
\frac{s_{i+1}\!-\!s_ip_i}{1\!+\!p_i},
\qquad  L_3(n):=\!s_1^2\!+\!\sum_{i=1}^{n-1} 
\frac{(s_{i+1}\!-\!s_ip_i)^2}{1\!-\!p_i^2},
\end{equation}
where \ $p_i:=\exp(-\beta d_i)$ \  and \ $d_i:=s_{i+1}-s_i, \ i=1,2, \ldots
,n-1$. \ To simplify calculations, we assume
that the first design point is at the origin, that is \ $s_1=0$, \ which does not change the general character of the presented results. Hence, in order to obtain the D-optimal design, one has to find
the maximum in \ $\boldsymbol d=(d_1,d_2,\ldots ,d_{n-1})$ \ of
\begin{equation}
  \label{Dopt}
{\mathcal D}(\boldsymbol d):=\det \big(\mathcal
I_{\alpha_0,\alpha_1}(n)\big)= L_1(n)L_3(n)-L_2^2(n),
\end{equation}
whereas K-optimal design minimizes the condition number \ ${\mathcal
  K}(\boldsymbol d)$  \ of \ $\mathcal I_{\alpha_0,\alpha_1}(n)$, \  where
\begin{equation}
  \label{Kopt}
{\mathcal K}(\boldsymbol d):=\frac 14 \bigg(L_1(n)+L_3(n)+\sqrt{\big(L_1(n)-L_3(n)\big)^2+4L_2^2(n)}\bigg)^2 \Big / \big(L_1(n)L_3(n)-L_2^2(n)\big). 
\end{equation}
Now, observe that 
\begin{equation*}
{\mathcal K}(\boldsymbol d)=g\big(\mathcal R(\boldsymbol d)\big), \qquad \text{where} \qquad g(x):=\frac 14 \big(\sqrt{x}+\sqrt{x-4}\big)^2, \quad x\geq 4,
\end{equation*} 
and
\begin{equation}
  \label{Ropt}
{\mathcal R}(\boldsymbol d):= \big(L_1(n)+L_3(n)\big)^2 \big / \big(L_1(n)L_3(n)-L_2^2(n)\big)\geq 4.
\end{equation}
As \ $g(x)$ \ is strictly monotone increasing, K-optimal design can be found by minimizing the objective function \ $\mathcal R(\boldsymbol d)$. \ Hence, the properties of K-optimal design for OU processes with linear trend are derived with the help of \ $\mathcal R(\boldsymbol d)$.

General results on D-optimal designs for models driven by OU processes
have already been formed and published \citep{ks,zba09}, but the
dependence of \ ${\mathcal R}(\boldsymbol d)$ \ on the design points
is far more complicated. Hence, in the next sections we investigate
some special cases in order to highlight the main differences between
the two design criteria. 

\begin{ex}
   \label{ex1}
Let the design space be \ ${\mathcal X}=[0,1]$ \ and consider a three-point restricted design \citep[see, e.g.,][]{bsschemo} where \ $s_1=0, \ s_2:=d, \ s_3=1$ \ with \ $0\leq d\leq 1$. \ In this case the objective functions \eqref{Dopt} and \eqref{Ropt} are univariate functions of \ $d$ \ and take the forms
\begin{equation*}
{\mathcal D}(d)=2\frac{(1-{\mathrm e}^{-\beta d})+d({\mathrm e}^{-\beta d}-{\mathrm e}^{-\beta (1-d)})-d(1-d)(1-{\mathrm e}^{-\beta})}{(1-{\mathrm e}^{-2\beta d})(1-{\mathrm e}^{-2\beta (1-d)})} \qquad \text{and} \qquad 
{\mathcal R}(d)= \frac {{\mathcal R}_1^2(d)}{{\mathcal R}_2(d)},
\end{equation*}
respectively, where
\begin{align*}
{\mathcal R}_1(d):=&\, (1-{\mathrm e}^{-\beta (1-d)})^2(1-{\mathrm e}^{-2\beta d})+(1-{\mathrm e}^{-\beta d})^2(1-{\mathrm e}^{-2\beta (1-d)})+(1-{\mathrm e}^{-2\beta d})(1-{\mathrm e}^{-2\beta (1-d)}) \\
&+d^2(1-{\mathrm e}^{-2\beta (1-d)})+(1-d{\mathrm e}^{-\beta (1-d)})^2(1-{\mathrm e}^{-2\beta d}), \\[2mm]
{\mathcal R}_2(d):= &\, 
2(1-{\mathrm e}^{-2\beta d})(1-{\mathrm e}^{-2\beta (1-d)})\big(1-{\mathrm e}^{-\beta d}+d^2(1-{\mathrm e}^{-\beta})-d(1+{\mathrm e}^{-\beta (1-d)})(1-{\mathrm e}^{-\beta d})\big).
\end{align*}
Direct calculations show that for all \ $\beta>0$ \ function \ $\mathcal D(d)$ \ has its maximum at \ $d=1/2$, \ that is the D-optimal three point restricted design is equidistant.

In case of K-optimality, the situation is completely different. Assume first \ $\beta^*\leq\beta\leq\beta ^{**}$, \ where \ $ \beta^*\approx 0.5718$ \ and \  $\beta^{**}\approx 4.9586$ \ are the only positive roots of \ $S(\beta)=0$, \ with
\begin{equation*}
S(\beta)\!:=\!(\beta^2-6\beta+4){\mathrm e}^{4\beta}\!+
(6\beta^2+6\beta-10){\mathrm e}^{3\beta}\!-(11\beta^2-10\beta-2)
{\mathrm e}^{2\beta}\!+(2\beta^2-6\beta+10){\mathrm e}^{\beta}-2\beta^2-4\beta-6.
\end{equation*}
Since
\begin{equation*}
\lim_{d\to 0}{\mathcal R}'(d)=-\frac {3{\mathrm e}^{\beta}-2}
{2\beta{\mathrm e}^{\beta}({\mathrm e}^{2\beta}-1)^2}S(\beta),
\end{equation*}
$\beta^*$ \ and \ $\beta^{**}$ \ are the only solutions of \ $\lim_{d\to 0}{\mathcal R}'(d)=0$, \ too. For \ $\beta \in[\beta^*,\beta^{**}]$, \ function  \ ${\mathcal R}(d)$ \ has a single extremal point in the \ $]0,1[$ \ interval, which corresponds to a maximum. Hence, as
\begin{equation*}
\lim_{d\to 0}{\mathcal R}(d)=\lim_{d\to 1}{\mathcal R}(d)=\frac {\big (3{\mathrm e}^{\beta}-2\big)^2}
{{\mathrm e}^{2\beta}-1},
\end{equation*}
the minimum of \ ${\mathcal R}(d)$ \ is reached at the boundary points $0$ and $1$, so the K-optimal design collapses. In contrast, for parameter values outside the interval \ $[\beta^*,\beta^{**}]$ \ K-optimal designs exist.
Figures \ref{fig1}a and \ref{fig1}b display the K-optimal value \ $d_{opt}$ \ plotted against the parameter \ $\beta$ \ for  intervals \ $]0,\beta^*[$ \ and  \ $]\beta^{**},100]$. \ We remark that \ $d_{opt}\to 0$ \ as \ $\beta \to \infty$, \ and the limit \ $0$ \ is the minimum point of 
\begin{equation*}
\lim_{\beta\to\infty}{\mathcal R}(d)=\frac {\big(d^4+4\big)^2}{2(d^2-d+1)}.
\end{equation*}

\begin{figure}[t]
\begin{center}
\leavevmode
\hbox{
\epsfig{file=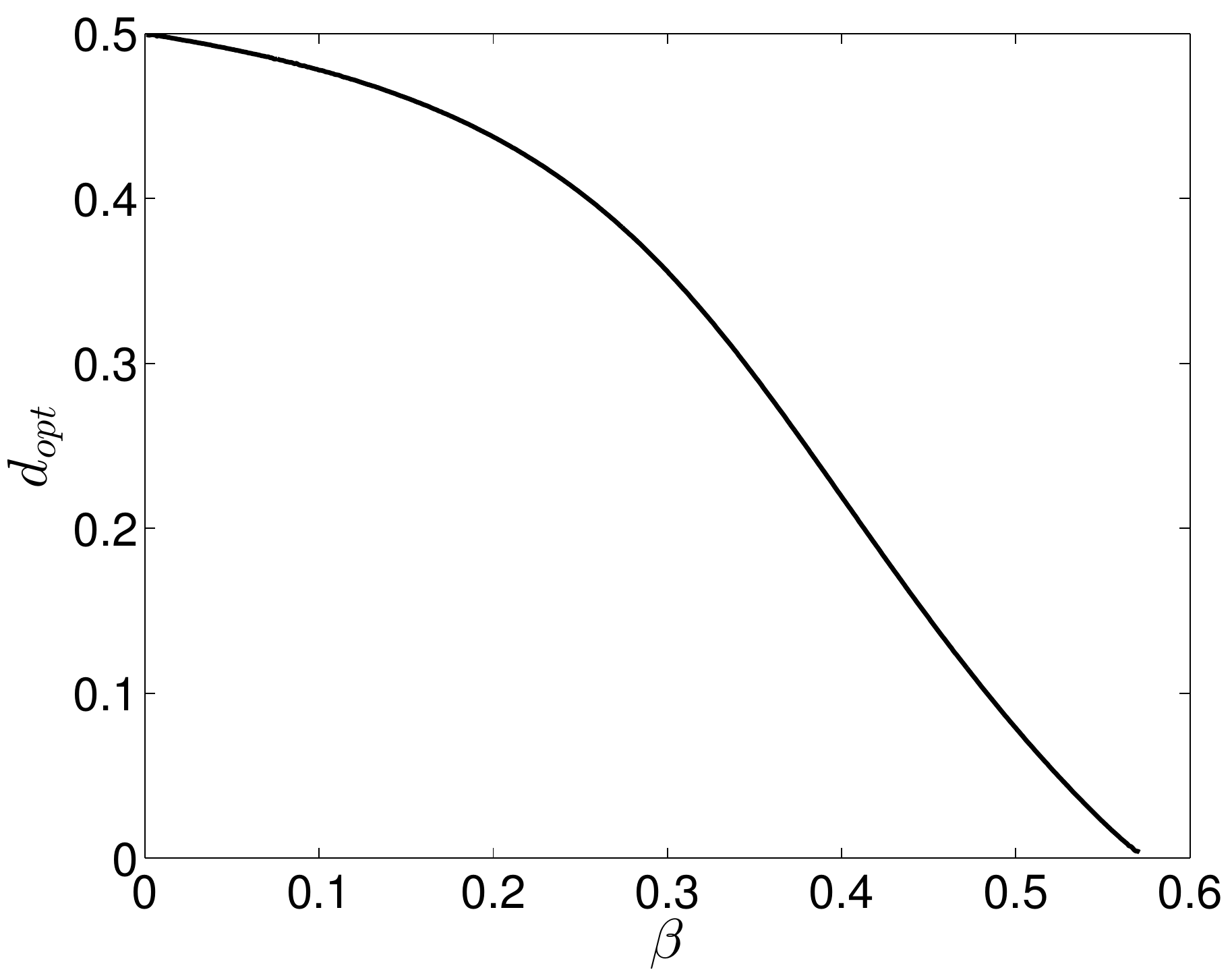,height=6cm, width=8 cm} \quad
\epsfig{file=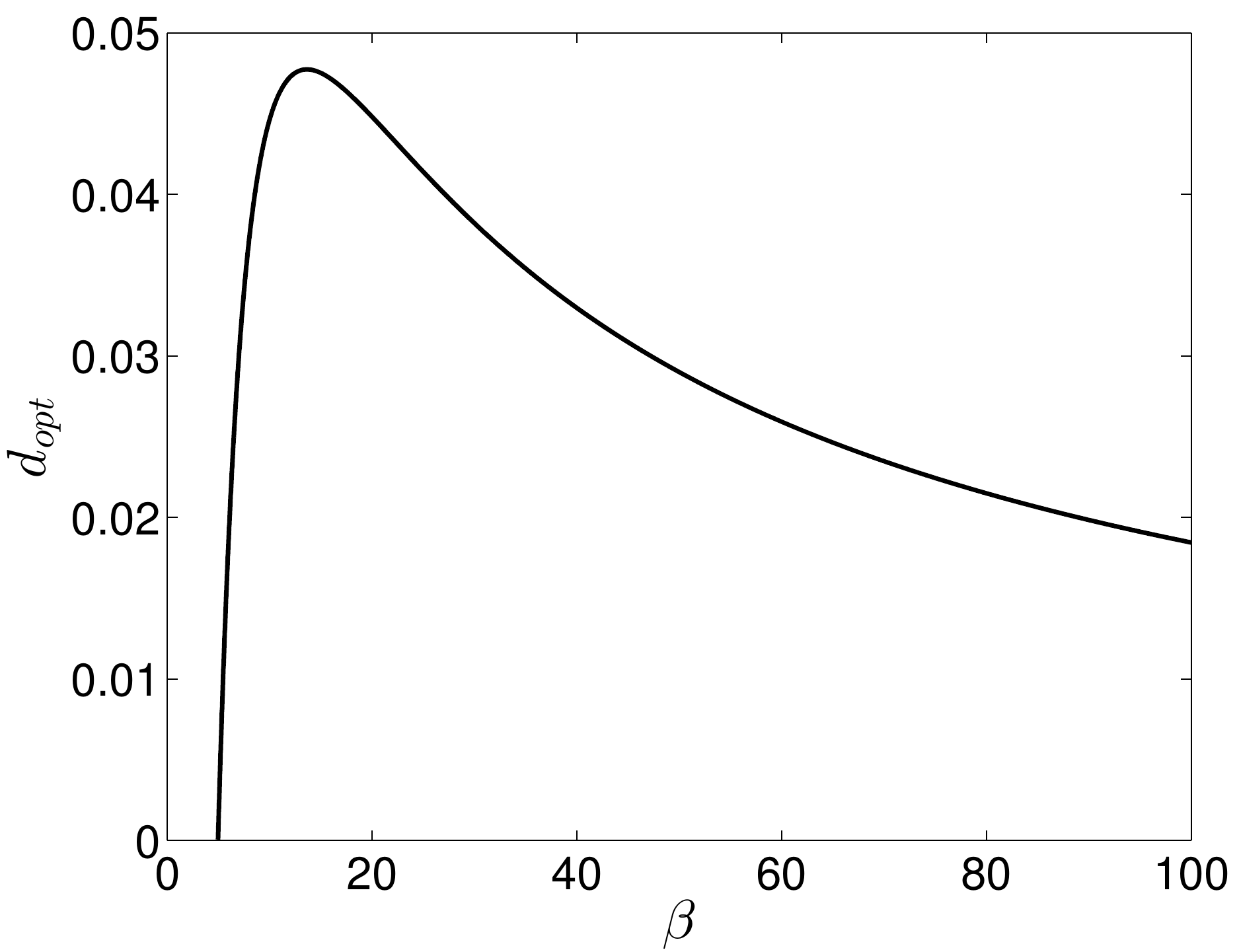,height=6cm, width=8 cm}}

\centerline{\hbox to 9 truecm {\scriptsize (a) \hfill (b)}}

\end{center}
\caption{K-optimal value \ $d_{opt}$ \ for the three point design \ $\boldsymbol\xi=\{0,d,1\}$ \ plotted against the parameter \ $\beta$ \ for the  intervals (a)  \ $]0,\beta^*[, \ \beta^*\approx 0.5718$, \ and  (b) \ $]\beta^{**},100], \ \beta^{**}\approx 4.9586$.}
\label{fig1}
\end{figure}

\end{ex}

\subsection{Optimality of increasing domain equidistant designs}
  \label{subs:subs2.1}

Consider an equidistant increasing domain design with step size \ $d>0$, \
that is the observation points are \ $\boldsymbol \xi=\big\{0,d,2d, \ldots
,(n-1)d\big\}$. In this case, \ $p_i=\exp (-\beta d), \ i=1,2, \ldots
,n-1,$ \ so the expressions in \eqref{entry} reduce to
\begin{align}
   \label{entry_d}
L_1(n)&=\frac {2-n+n{\mathrm e}^{\beta d}}{{\mathrm e}^{\beta
  d}+1}, \qquad \qquad L_2(n)=\frac{d(n-1)}2 L_1(n), \\
L_3(n)&=\frac {d^2(n-1)}{{\mathrm e}^{2\beta
  d}-1}\bigg(\frac {n(2n-1)({\mathrm e}^{\beta d}-1)^2}6+n\big({\mathrm e}^{\beta
  d}-1\big)+1 \bigg), \nonumber
\end{align} 
and the objective functions \ $\mathcal D$, \  $\mathcal K$ \ and \ $\mathcal R$ \ defined by  \eqref{Dopt}, \eqref{Kopt} and \eqref{Ropt}, respectively, are univariate functions of  \ $d$. 

\begin{thm}
  \label{equiresgen}
For model \eqref{model} with covariance structure
\eqref{oucov} and equidistant increasing domain design with step size \
$d>0$,  \ function  \ $\mathcal D$ \ is monotone
increasing in \ $d$, \ whereas \ $\mathcal R$ \ (and \ $\mathcal K$) \ has at least one global minimum point, that is, there exists a K-optimal design.
\end{thm}

\begin{figure}[t]
\begin{center}
\leavevmode
\hbox{
\epsfig{file=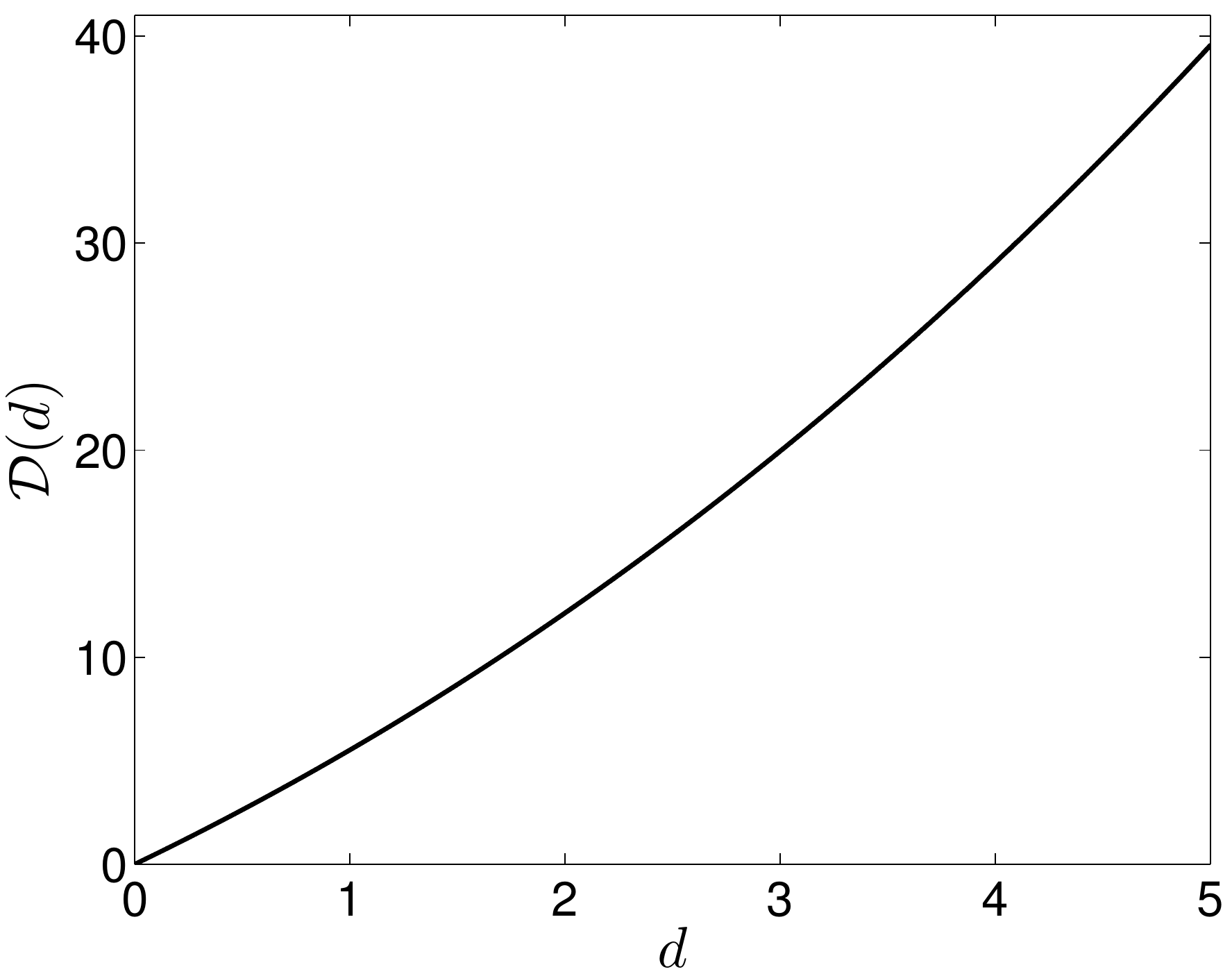,height=6cm, width=8 cm} \quad
\epsfig{file=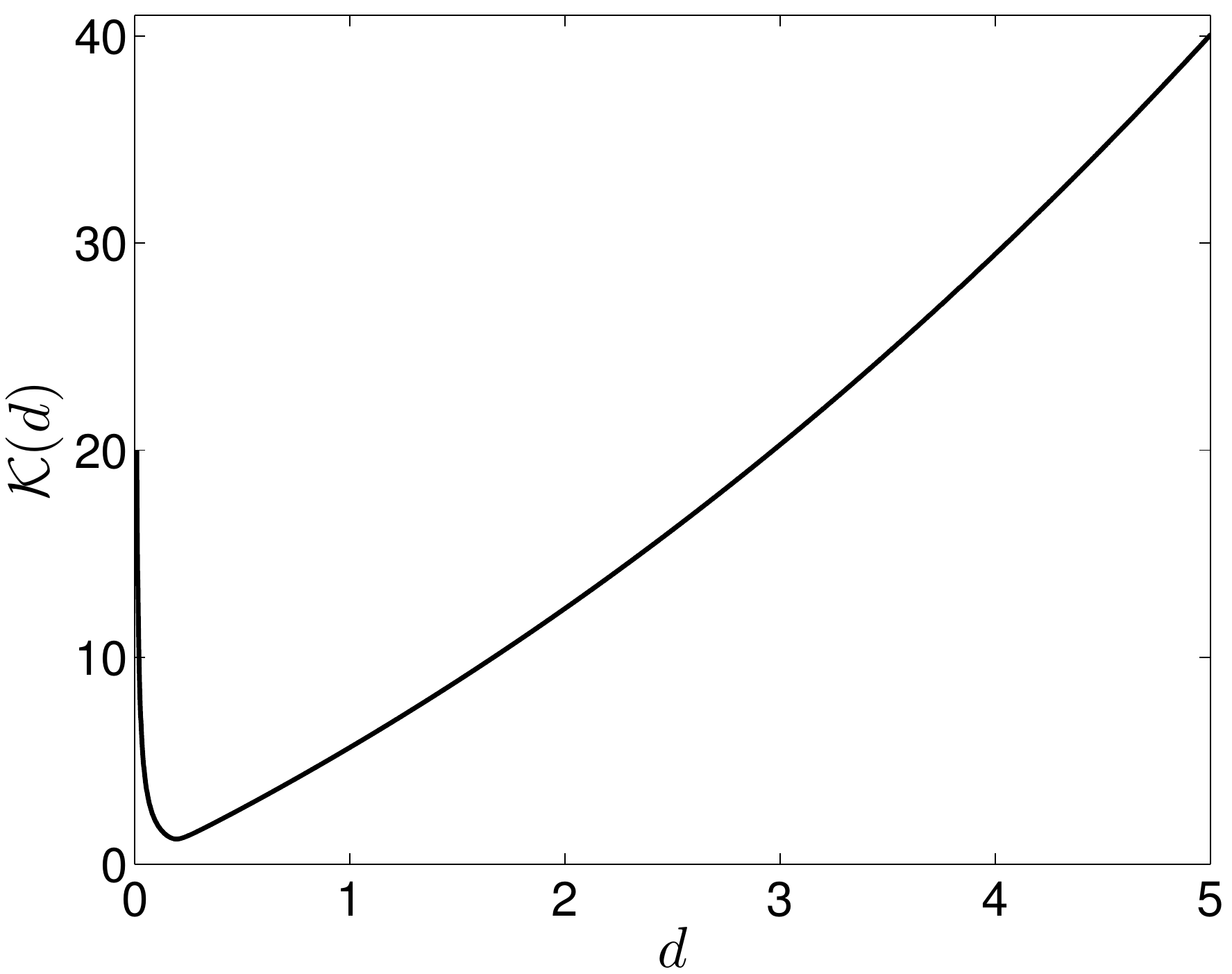,height=6cm, width=8 cm}}

\centerline{\hbox to 9 truecm {\scriptsize (a) \hfill (b)}}

\end{center}
\caption{Objective functions for two-point  D-optimal (a) and 
  K-optimal (b) designs with \ $\beta=0.1$.}
\label{fig2}
\end{figure}

As a special case consider a two-point design \ $\{0,d\}$. \ Figures
\ref{fig2}a and \ref{fig2}b show the behaviour of \ ${\mathcal D}(d)$
\ and \ ${\mathcal K}(d)$, \ respectively, for \ $\beta=0.1$, \
whereas the following theorem formulates a general result on the
two-point K-optimal design. 
\begin{thm}
   \label{eqirestwop}
For model \eqref{model} with covariance structure \eqref{oucov}, there
exists a unique K-optimal two-point design \ $\{0,d_{opt}\}$, \ where \
$d_{opt}$ \ is the unique solution of
\begin{equation}
  \label{twopeq}
\big(d^2-2\big){\mathrm e}^{3\beta d}+2(\beta d+1){\mathrm e}^{2\beta
  d}-\big(\beta d^3+d^2+2\beta d-2\big){\mathrm 
  e}^{\beta d}-2=0.
\end{equation}
\end{thm}

\subsection{Comparison of equidistant designs}
  \label{subs:subs2.2}
According to the ideas of \citet{hoel58}, we investigate the change
of D- and K-optimality criteria arising from doubling the number of
sub-intervals in an equidistant partition of a fixed design interval, and we also study the situation when the length of the design interval is also doubled. 
The former approach refers to infill-, whereas the latter to increasing domain asymptotics.

Let the design space be  \ ${\mathcal X}=[a,b]\subset {\mathbb R}, \ a<b$, \ and denote by \ $\widetilde{\mathcal X}$ \ the interval \ $[2a,2b]$. \ 
Obviously, without loss of generality we may assume  \
${\mathcal X}=[0,1]$ \ and consider the sequences \ $\boldsymbol \xi_n$ \ and \ 
$\boldsymbol \xi_{2n}$ \ of designs on \ ${\mathcal X}$, \ where \ $\boldsymbol \xi_n:=\big \{0,1/n, \ldots ,(n-1)/n,1\big\},\ n\geq 2, \ n\in{\mathbb N}$, 
\ and \ designs  \ $\widetilde{\boldsymbol \xi}_{2n}:=\big \{0,1/n, \ldots ,(2n-1)/n,2\big\}$ \ on \ $\widetilde{\mathcal X}=[0,2]$.
\begin{thm}
  \label{doublim}
For model \eqref{model} with covariance structure \eqref{oucov}
\begin{equation}
  \label{k_lim}
\lim_{n\to\infty}\frac {{\mathcal D}(2n)}{{\mathcal D}(n)}=
\lim_{n\to\infty}\frac {{\mathcal K}(2n)}{{\mathcal K}(n)}=1,  
\end{equation}
where \ $\mathcal D(n)$ \ and \ $\mathcal K(n)$ \ are the values of
the objective functions \eqref{Dopt} and \eqref{Kopt}, respectively,
corresponding to the design \ $\boldsymbol \xi_n$ \ on \ ${\mathcal X}$.  
\end{thm} 

Limits \eqref{k_lim} show that for both investigated design criteria, if one has a dense enough equidistant partition of a fixed design space, there is no use of doubling the number of intervals in the partition, which is in accordance with the results of \citet{hoel58} for classical polynomial regression and \citet{ks} for OU processes with a constant trend.

\begin{thm}
  \label{doubintlim}
For model \eqref{model} with covariance structure \eqref{oucov}
\begin{equation}
  \label{kint_lim}
\lim_{n\to\infty}\frac {\widetilde{\mathcal D}(2n)}{{\mathcal D}(n)}=D(\beta)
\qquad \text{and} \qquad
\lim_{n\to\infty}\frac {\widetilde{\mathcal K}(2n)}{{\mathcal K}(n)}=K(\beta),  
\end{equation}
with
\begin{align}
  \label{DKform}
D(\beta)&:= \frac{16(\beta+1)(\beta^2+3\beta +3)}{(\beta+2) (\beta^2+6\beta +12)} \qquad \text{and} \\
K(\beta)&:=\frac{(\beta+2)(\beta^2+6\beta +12) \big(7\beta^2+9\beta +3+\sqrt{37\beta^4+78\beta^3+51\beta^2+18\beta+9}\big)^2}{4(\beta+1) (\beta^2+3\beta +3)\big(4\beta^2+9\beta +3+\sqrt{13\beta^4+48\beta^3+33\beta^2-18\beta+9}\big)^2}, \nonumber
\end{align}
where \ $\mathcal D(n)$ \ and \ $\widetilde{\mathcal D}(2n)$ \ are the values of
the objective function \eqref{Dopt}, whereas  \ $\mathcal K(n)$ \ and \ $\widetilde{\mathcal K}(2n)$ \ are the values of
the objective function \eqref{Kopt} corresponding to the designs \ $\boldsymbol \xi_n$ \ and \ $\widetilde{\boldsymbol \xi}_{2n}$, \ respectively.
\end{thm} 

Note that \ $D(\beta)$ \ is strictly increasing with \ $\lim_{\beta\to 0}D(\beta)=2$ \ and \ $\lim_{\beta\to\infty}D(\beta)=16$, \ whereas \ $\lim_{\beta\to 0}K(\beta)=2$, \ $K(\beta)$ \ has a single maximum of $2.3454$ at $0.2730$, and then it is strictly decreasing with \ $\lim_{\beta\to\infty}K(\beta)=\big(7+\sqrt{37}\big)^2/\big(8+2\sqrt{13}\big)^2\approx 0.7397$. \ Hence, doubling the interval over which the dense enough equidistant observations are made at least doubles the information on the unknown regression parameters \ $(\alpha_0,\alpha_1)$, \ which supports the extension of the design space. Moreover, after the maximum point of \ $K(\beta)$ \ the larger the covariance parameter \ $\beta$, \ the more we gain in efficiency in terms of the condition number with extending the interval where the observations are made.

\section{Ornstein-Uhlenbeck sheets with linear trend}
   \label{sec:sec3}

As a spatial generalization of model \eqref{model}, consider the spatial process
\begin{equation}
   \label{2Dmodel}
Y(s,t)=\alpha_0+\alpha_1s+\alpha_2t+U(s,t),
\end{equation}
where the design points are taken from a compact design space
\ $\mathcal{X}=[a_1,b_1]\times [a_2, b_2]$, \ with \ $b_1>a_1$ \ and \
$b_2>a_2$, \ and \ $U (s,t), \ s,t\in {\mathbb R}$, \ is a
stationary OU sheet, i.e.,  
a centered Gaussian process with covariance structure
\begin{equation}
   \label{2Doucov}
{\mathsf E}\,U(s_1,t_1)U(s_2,t_2)=
\frac{\sigma^2}{4\beta\gamma}\exp\big
(-\beta|s_1-s_2|-\gamma|t_1-t_2|\big ),
\end{equation}
where \ $\beta>0, \ \gamma>0, \ \sigma>0$. \ Similar to the OU process, \
$U(s,t)$ \ can be represented as
\begin{equation}
   \label{2Dourep}
U(s,t)=\frac{\sigma}{2\sqrt{\beta\gamma}}{\mathrm
  e}^{-\beta s-\gamma t}{\mathcal W}\big({\mathrm e}^{2\beta s},
  {\mathrm e}^{2\gamma t}\big),
\end{equation}
where \ ${\mathcal W}(s,t), \ s,t\in {\mathbb R}$, \ is a  standard
  Brownian sheet \citep{bpz, bs}. Again, we assume that the parameters \ $\beta, \ \gamma$ \ and \ $\sigma$ \ of the OU sheet \ $U(s,t)$ \ driving model \eqref{2Dmodel} are known.

We investigate regular grid designs of the form \ $\big\{(s_i,t_j): \ i=1,2, \ldots ,n, \ j=1,2, \ldots ,m\big\}\subset \mathcal{X}=[a_1,b_1]\times [a_2,
b_2], \ n,m\geq 2$, \ and without loss of generality we may assume \
$a_1\leq s_1<s_2< \ldots < s_n\leq b_1$ \ and \ $a_2\leq t_1<t_2< \ldots < t_m\leq b_2$ \ \citep{bsschemo}, and that \ $U(s,t)$ \ has a unit variance. Again, the general form of the FIM \ $\mathcal I_{\alpha_0,\alpha_1,\alpha_2}(n,m)$ \ on regression parameters \ $\alpha_0,\ \alpha_1$ \ and \ $\alpha_2$ \ of model \eqref{2Dmodel} based on observations \ $\big\{Y(s_i,t_j), \ i=1,2,\ldots, n, \ j=1,2,\ldots, m\big\}$ \ equals
\begin{equation*}
\mathcal I_{\alpha_0,\alpha_1,\alpha_2}(n,m)=G(n,m)C(n,m)^{-1}G(n,m)^{\top},
\end{equation*} 
where \ $C(n,m)$ \ denotes the covariance matrix of the observations and
\begin{equation*}
G(n,m):=\begin{bmatrix} 
1 & 1 & \cdots & 1 & 1 & 1 & \cdots & 1 & \cdots & 1 & 1 & \cdots & 1 \\
s_1 & s_1 & \cdots & s_1 & s_2 & s_2 & \cdots & s_2 & \cdots & s_n & s_n & \cdots & s_n \\
t_1 & t_2 & \cdots & t_m & t_1 & t_2 & \cdots & t_m & \cdots & t_1 & t_2 &\cdots & t_m \\
\end{bmatrix}.
\end{equation*}
The following theorem gives the exact form of the FIM \ $I_{\alpha_0,\alpha_1,\alpha_2}(n,m)$. 

\begin{thm}
  \label{FIMsheet}
Consider the OU model \eqref{2Dmodel} with covariance structure \eqref{2Doucov} observed in points \ $\big\{(s_i,t_j): \ i=1,2, \ldots ,n, \ j=1,2, \ldots ,m\big\}$. \ Then
\begin{equation}
  \label{FIMform}
\mathcal I_{\alpha_0,\alpha_1,\alpha_2}(n,m) =\begin{bmatrix}
L_1(n)M_1(m) & L_2(n)M_1(m) & L_1(n)M_2(m) \\
L_2(n)M_1(m) & L_3(n)M_1(m) & L_2(n)M_2(m) \\
L_1(n)M_2(m) & L_2(n)M_2(m) & L_1(n)M_3(m) 
\end{bmatrix}
\end{equation}
with
\begin{alignat*}{3}
L_1(n)&:=\!1\!+\!\sum_{i=1}^{n-1} \frac{1\!-\!p_i}{1\!+\!p_i}, \qquad \
L_2(n)&&:=\!s_1\!+\!\sum_{i=1}^{n-1} 
\frac{s_{i+1}\!-\!s_ip_i}{1\!+\!p_i},
\qquad  L_3(n)&&:=\!s_1^2\!+\!\sum_{i=1}^{n-1} 
\frac{(s_{i+1}\!-\!s_ip_i)^2}{1\!-\!p_i^2},\\
M_1(m)&:=\!1\!+\!\sum_{i=1}^{m-1} \frac{1\!-\!q_i}{1\!+\!q_i}, \qquad
M_2(m)&&:=\!t_1\!+\!\sum_{i=1}^{m-1} 
\frac{t_{i+1}\!-\!t_iq_i}{1\!+\!q_i},
\qquad M_3(m)&&:=\!t_1^2\!+\!\sum_{i=1}^{m-1} 
\frac{(t_{i+1}\!-\!t_iq_i)^2}{1\!-\!q_i^2},
\end{alignat*}
where \ $p_i:=\exp(-\beta d_i)$ \ with  \  $d_i:=s_{i+1}-s_i, \ \ i=1,2, \ldots ,n-1$, \ and \ $q_j:=\exp(-\gamma \delta_j)$ \ with 
$\delta_j:=t_{j+1}-t_j, \ j=1,2, \ldots, m-1$.
\end{thm}

Again, to simplify calculations we assume \ $s_1=t_1=0$, \ so the
D-optimal design maximizes 
\begin{equation}
  \label{2DDopt}
{\mathcal D}(\boldsymbol d,\boldsymbol \delta):=\det \big(\mathcal I_{\alpha_0,\alpha_1,\alpha_2}(n,m)\big)=
L_1(n)M_1(m)\big(L_1(n)L_3(n)-L_2^2(n)\big)\big(M_1(m)M_3(m)-M_2^2(m)\big)
\end{equation}
both in \ $\boldsymbol d=(d_1,d_2,\ldots ,d_{n-1})$ \ and \ $\boldsymbol \delta=(\delta_1,\delta_2,\ldots ,\delta_{m-1})$, \ whereas to obtain the K-optimal design one has to minimize the condition number \ ${\mathcal K}(\boldsymbol d,\boldsymbol \delta)$ \ of \ $\mathcal I_{\alpha_0,\alpha_1,\alpha_2}(n,m)$. \ Using the expressions of \citet{smith} for the eigenvalues of a $3\times 3$ symmetric matrix, one can easily show
\begin{equation}
  \label{2DKopt}
{\mathcal K}(\boldsymbol d,\boldsymbol \delta)=
\frac{\tr \big(\mathcal I_{\alpha_0,\alpha_1,\alpha_2}(n,m)\big)+ 
\sqrt{6\tr \big(\mathcal I_{\alpha_0,\alpha_1,\alpha_2}^{\,2}(n,m)\big) 
-2\tr^2 \big(\mathcal I_{\alpha_0,\alpha_1,\alpha_2}(n,m)\big)}\cos (\varphi)}
{\tr \big(\mathcal I_{\alpha_0,\alpha_1,\alpha_2}(n,m)\big)+ 
\sqrt{6\tr \big(\mathcal I_{\alpha_0,\alpha_1,\alpha_2}^{\,2}(n,m)\big) 
-2\tr^2 \big(\mathcal I_{\alpha_0,\alpha_1,\alpha_2}(n,m)\big)}\cos \big(\varphi+2\pi /3 \big)},
\end{equation}
where \ $\varphi := \frac 13 \arccos (\varrho ) \in [0,\pi/3]$, \ with
\begin{equation*}
  \varrho:=\frac{54\det \big(\mathcal I_{\alpha_0,\alpha_1,\alpha_2}(n,m)\big)
+\tr\big(\mathcal I_{\alpha_0,\alpha_1,\alpha_2}(n,m)\big) \Big(
9\tr \big(\mathcal I_{\alpha_0,\alpha_1,\alpha_2}^{\,2}(n,m)\big) 
-5\tr^2 \big(\mathcal I_{\alpha_0,\alpha_1,\alpha_2}(n,m)\big)\Big)}
{\sqrt 2  \Big(
3\tr \big(\mathcal I_{\alpha_0,\alpha_1,\alpha_2}^{\,2}(n,m)\big) 
-\tr^2 \big(\mathcal I_{\alpha_0,\alpha_1,\alpha_2}(n,m)\big)\Big)^{3/2}}.
\end{equation*}

\begin{ex}
   \label{ex2}
Let the design space be the unit square \ $\mathcal X=[0,1]^2$ \ and consider a nine-point restricted regular grid design, where \ $s_1=t_1=0, \ s_2=d, \ t_2=\delta, \ s_3=t_3=1$, \ with \ $0\leq d,\delta \leq 1$. \ In this case
the  FIM \eqref{FIMsheet} equals
\begin{equation*}
\mathcal I_{\alpha_0,\alpha_1,\alpha_2}(d,\delta) =\begin{bmatrix}
L_1(d)M_1(\delta ) & L_2(d)M_1(\delta ) & L_1(d)M_2(\delta ) \\
L_2(d)M_1(\delta ) & L_3(d)M_1(\delta ) & L_2(d)M_2(\delta ) \\
L_1(d)M_2(\delta ) & L_2(d)M_2(\delta ) & L_1(d)M_3(\delta ) 
\end{bmatrix}
\end{equation*}
with 
\begin{alignat*}{3}
L_1(d)&\!:=\!\frac 2{1\!+\!{\mathrm e}^{-\beta d}}\!+\!\frac {1\!-\!{\mathrm e}^{-\beta (1-d)}}{1\!+\!{\mathrm e}^{-\beta (1-d)}},  \ 
L_2(d)&&\!:=\!\frac d{1\!+\!{\mathrm e}^{-\beta d}}\!+\!\frac {1\!-\!d{\mathrm e}^{-\beta (1-d)}}{1\!+\!{\mathrm e}^{-\beta (1-d)}},  \ 
L_3(d)&&\!:=\!\frac {d^2}{1\!-\!{\mathrm e}^{-2\beta d}}\!+\!\frac {(1\!-\!d{\mathrm e}^{-\beta (1-d)})^2}{1\!-\!{\mathrm e}^{-2\beta (1-d)}}, \\ 
M_1(\delta)&\!:=\!\frac 2{1\!+\!{\mathrm e}^{-\gamma \delta}}\!+\!\frac {1\!-\!{\mathrm e }^{-\gamma (1-\delta)}}{1\!+\!{\mathrm e}^{-\gamma (1-\delta)}}, \ 
M_2(\delta)&&\!:=\!\frac {\delta}{1\!+\!{\mathrm e}^{-\gamma \delta}}\!+\!\frac {1\!-\!\delta {\mathrm e }^{-\gamma (1-\delta)}}{1\!+\!{\mathrm e}^{-\gamma (1-\delta)}}, \ 
M_3(\delta)&&\!:=\!\frac {\delta^2}{1\!-\!{\mathrm e}^{-2\gamma \delta}}\!+\!\frac {(1\!-\!\delta {\mathrm e }^{-\gamma (1-\delta)})^2}{1\!-\!{\mathrm e}^{-2\gamma (1-\delta)}},
\end{alignat*}
so both the determinant \ $\mathcal D$ \ and the condition number \ $\mathcal K$ \ of \ $\mathcal I_{\alpha_0,\alpha_1,\alpha_2}(d,\delta)$ \  are bivariate functions of \ $d$ \ and \ $\delta$. \ As for all possible parameter values \ $\beta>0$ \ function \ $L_1(d)$ \ reaches its unique maximum at \ $d=1/2$, \ and obviously the same holds for \ $M_1(\delta)$ \ for all \ $\gamma>0$, \ representations \eqref{Dopt} and \eqref{2DDopt} together with the results of Example \ref{ex1} imply that the D-optimal nine-point restricted regular grid design is directionally equidistant.   

Similar to the temporal case of Example \ref{ex1}, a non-collapsing K-optimal design exists only outside a certain region of the \ $(\beta,\gamma)$ \ parameter space. In Figures \ref{fig3}a and \ref{fig3}b  the \ $d_{opt}$ \ and \ $\delta_{opt}$ \ coordinates of the minimum point of \ $\mathcal K(d,\delta)$ are plotted as functions of parameters, where $0$ values correspond to collapsing designs, whereas Figures \ref{fig3}c and \ref{fig3}d display the corresponding contour plots.

\begin{figure}[t]
\begin{center}
\leavevmode
\hbox{
\epsfig{file=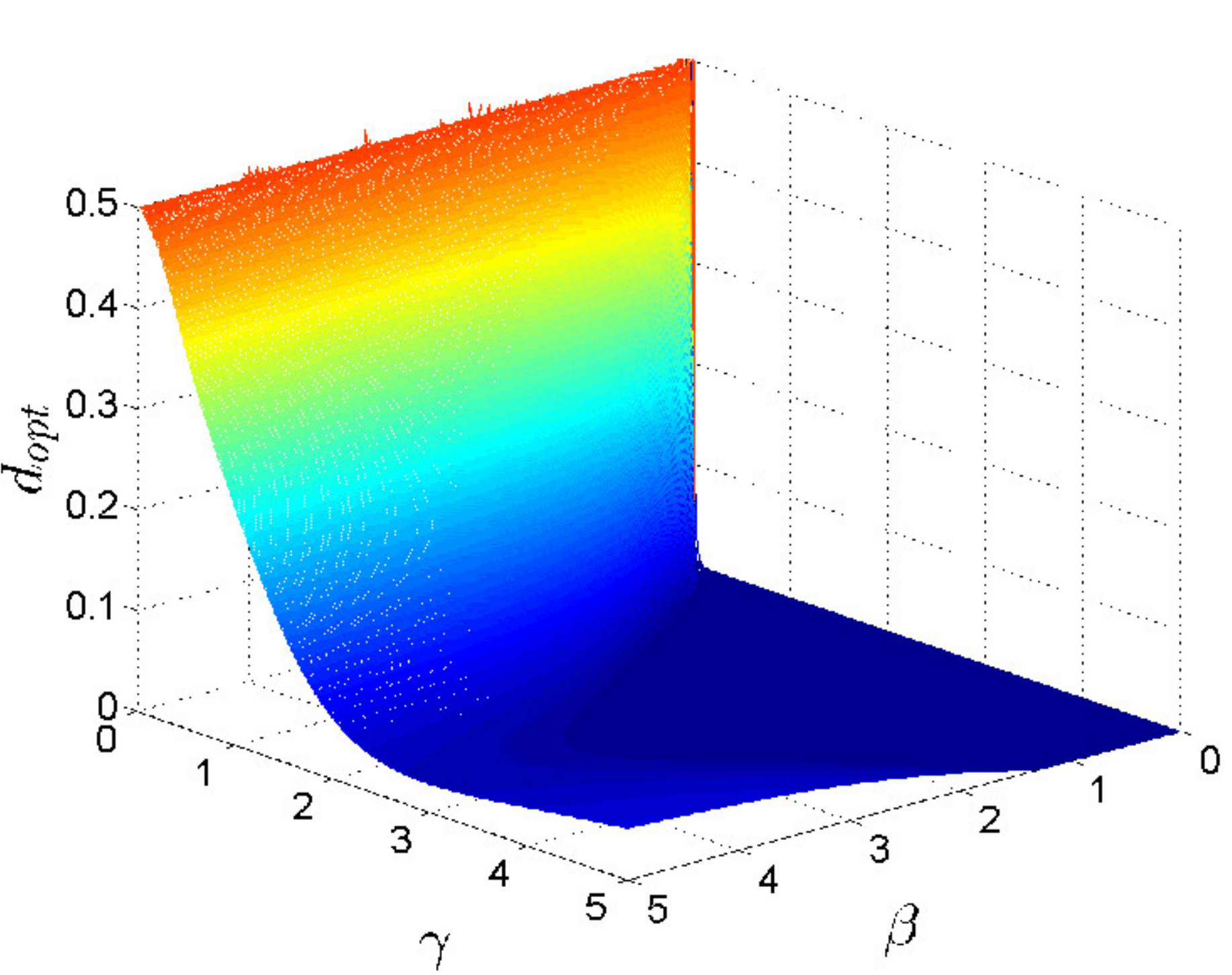,height=6cm, width=8 cm} \quad
\epsfig{file=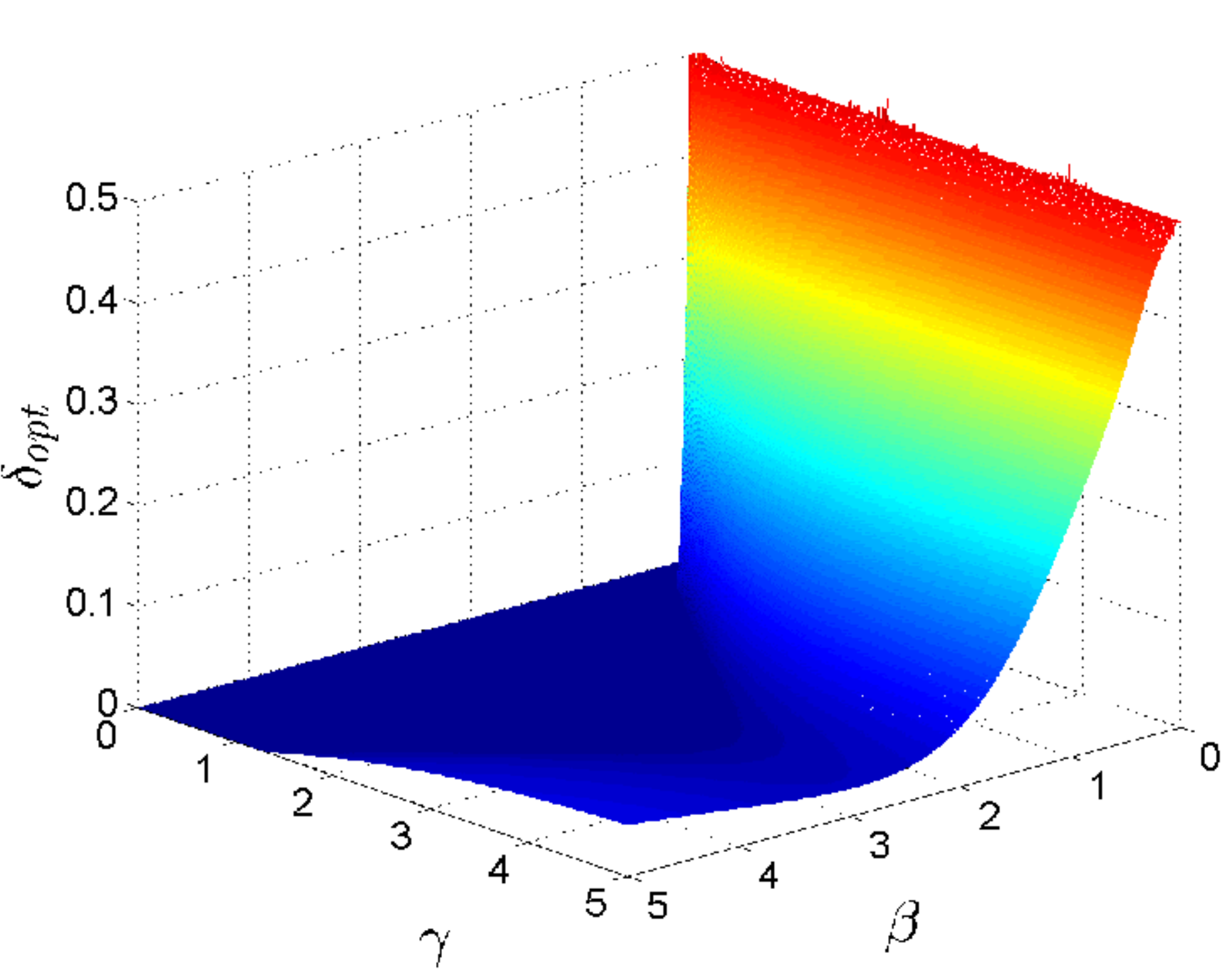,height=6cm, width=8 cm}}

\centerline{\hbox to 9 truecm {\scriptsize (a) \hfill (b)}}

\smallskip

\leavevmode
\hbox{
\epsfig{file=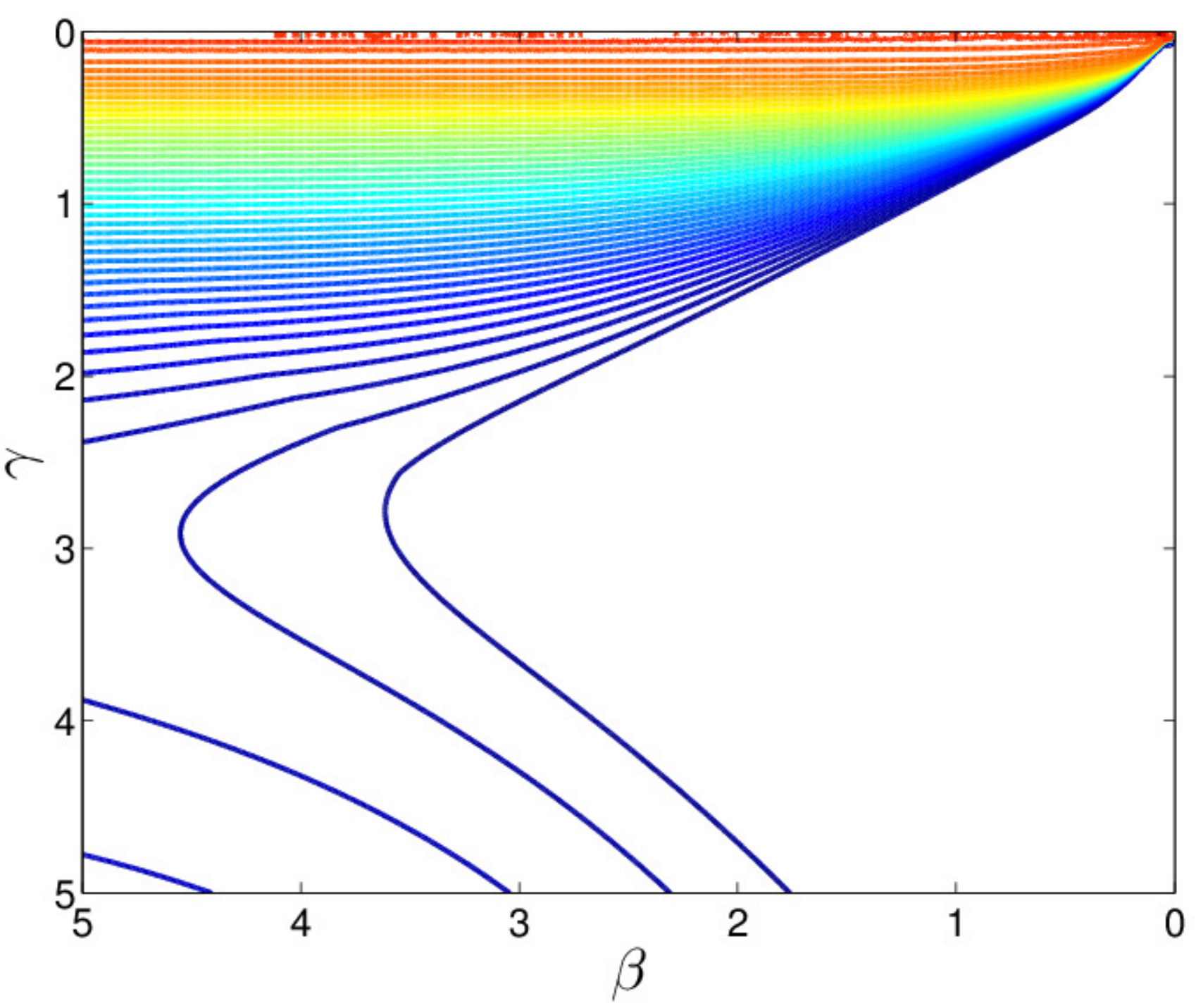,height=6cm, width=8 cm} \quad
\epsfig{file=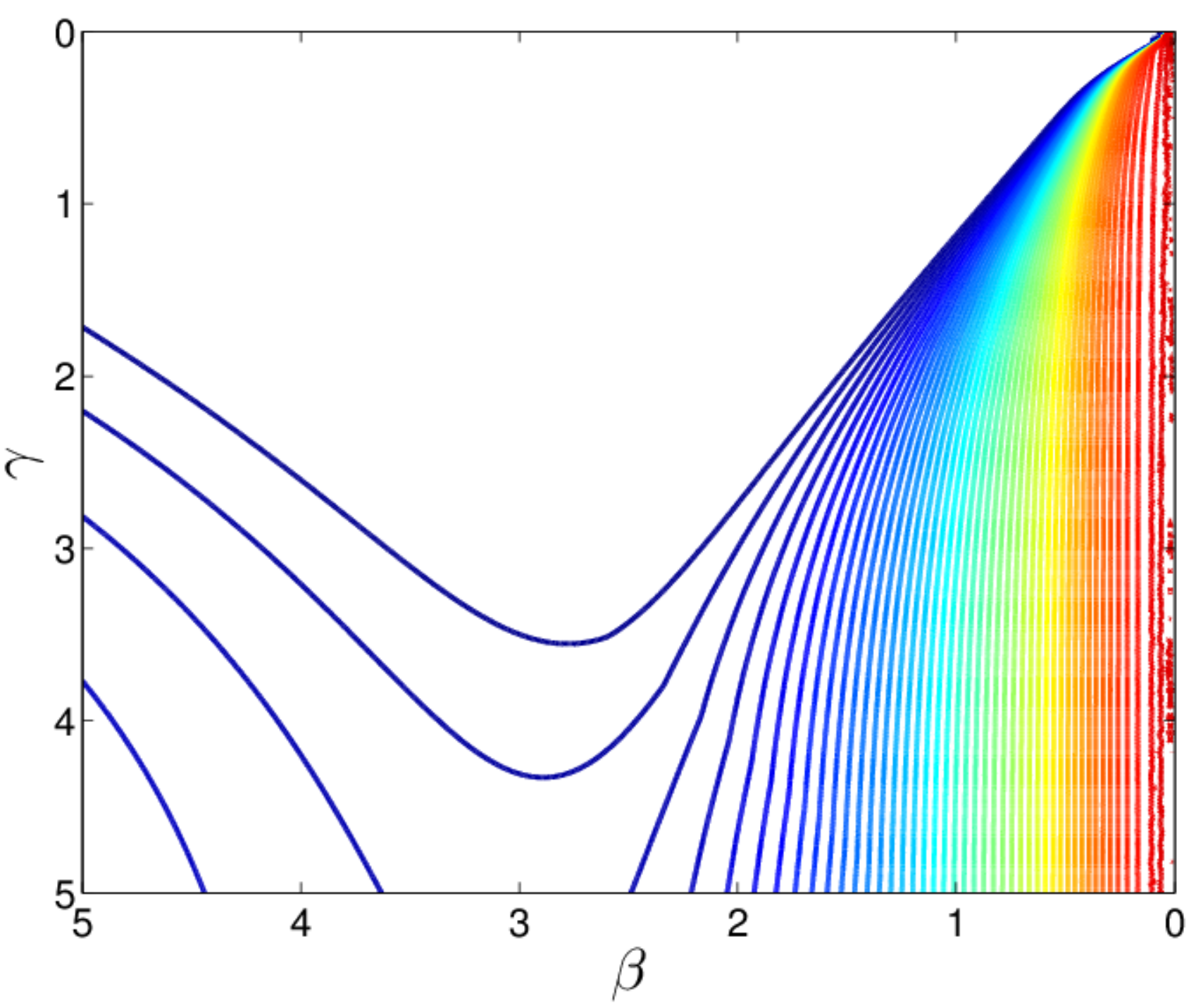,height=6cm, width=8 cm}}

\centerline{\hbox to 9 truecm {\scriptsize (c) \hfill (d)}}

\end{center}
\caption{ $d_{opt}$  (a) and \ $\delta_{opt}$  (b) coordinates of the minimum point of the objective function \ $\mathcal K(d,\delta)$ \ and the corresponding contour plots \ ($d_{opt}$:  (c); \ $\delta_{opt}$:  (d)) for a nine-point restricted regular grid design.}
\label{fig3}
\end{figure}
\end{ex}

\subsection{Optimality of increasing domain equidistant designs}
  \label{subs:subs3.1}

Consider now the directionally equidistant regular grid design with step sizes \ $d>0$ \ and \ $\delta>0$ \ consisting of observation locations  \ $\boldsymbol \xi = \big\{\big((i\!-\!1)d,(j\!-\!1)\delta\big) : \ i=1,2, \ldots, n, \ j=1,2,\ldots ,m\big\}$. \ In this situation, \ $L_1(n), \ L_2(n)$ \ and \ $L_3(n)$ have forms given by \eqref{entry_d}, and in a similar way we have
\begin{align}
   \label{entry_m}
M_1(m)&=\frac {2-m+m{\mathrm e}^{\gamma\delta}}{{\mathrm e}^{\gamma
  \delta}+1}, \qquad \qquad M_2(m)=\frac{\delta(m-1)}2 M_1(m), \\
M_3(m)&=\frac {\delta^2(m-1)}{{\mathrm e}^{2\gamma
  \delta}-1}\bigg(\frac {m(2m-1)({\mathrm e}^{\gamma \delta}-1)^2}6+m\big({\mathrm e}^{\gamma \delta}-1\big)+1 \bigg). \nonumber
\end{align} 
Hence, objective functions \ $\mathcal D$ \ and \ $\mathcal K$ \ are bivariate
functions of  \ $d$ \ and \ $\delta$. 

\begin{thm}
  \label{2Dequiresgen}
For model \eqref{2Dmodel} with covariance structure \eqref{2Doucov} and directionally equidistant increasing domain design with step sizes \
$d>0$  \ and \ $\delta>0$, \ function  \ $\mathcal D$ \ is monotone
increasing both in \ $d$ \ and \ $\delta$.
\end{thm}

Similar to Section \ref{subs:subs2.1}, in case of K-optimality one faces a different situation.

\begin{ex}
  \label{ex3}
Consider the four-point increasing domain regular grid design \ $\big\{(0,0),(d,0),(0,\delta),(d,\delta)\big\}$ \ for the process \eqref{2Dmodel} with parameters \ $\beta=0.2, \ \gamma=0.3$. \ By Theorem \ref{2Dequiresgen}, there is no D-optimal design,  whereas Figures \ref{fig4}a and \ref{fig4}b showing the objective function \ $\mathcal K(d,\delta)$ \ and the corresponding contour plot, respectively, clearly indicate the existence of a K-optimal design.

\begin{figure}[t]
\begin{center}
\leavevmode
\hbox{
\epsfig{file=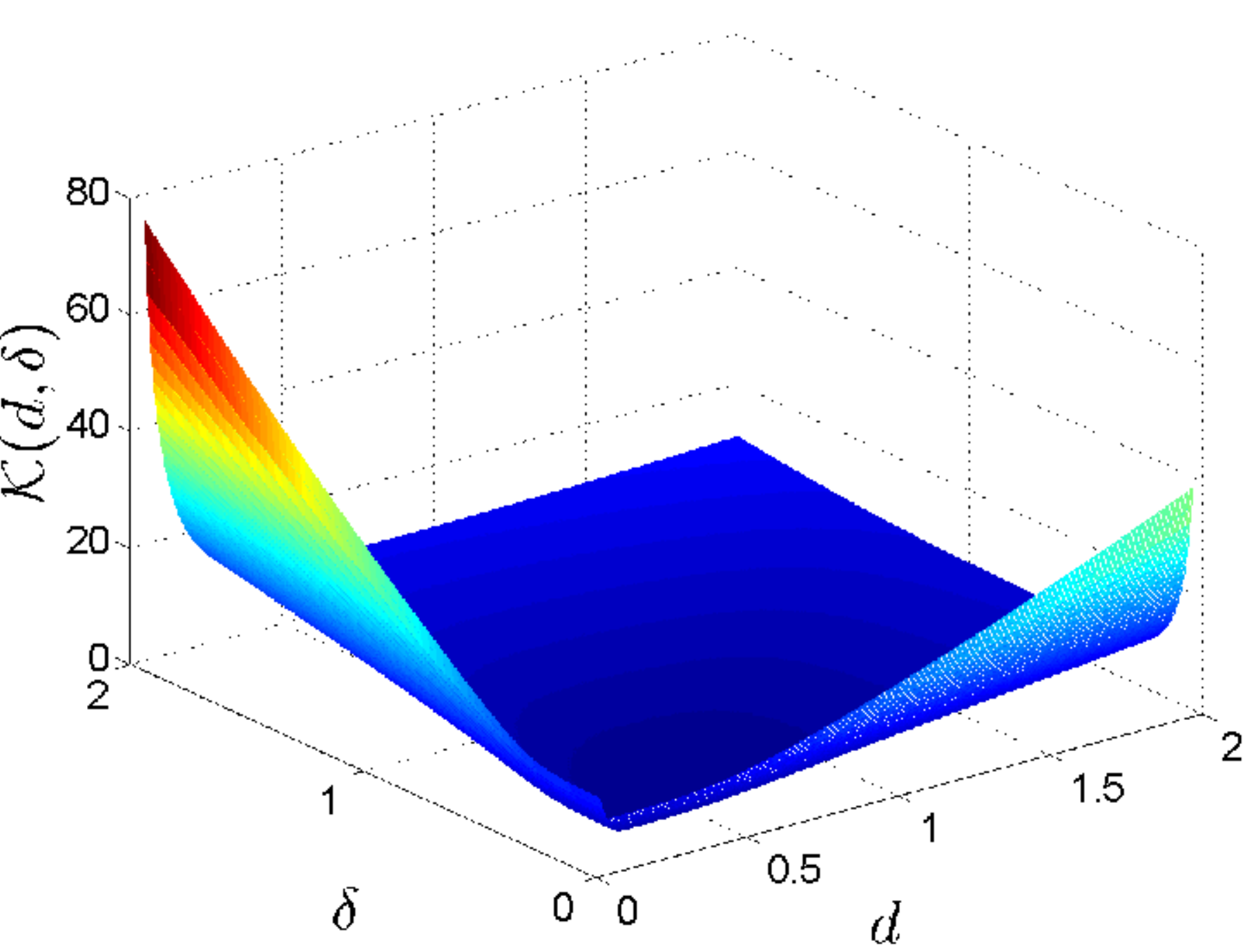,height=6cm, width=8 cm} \quad
\epsfig{file=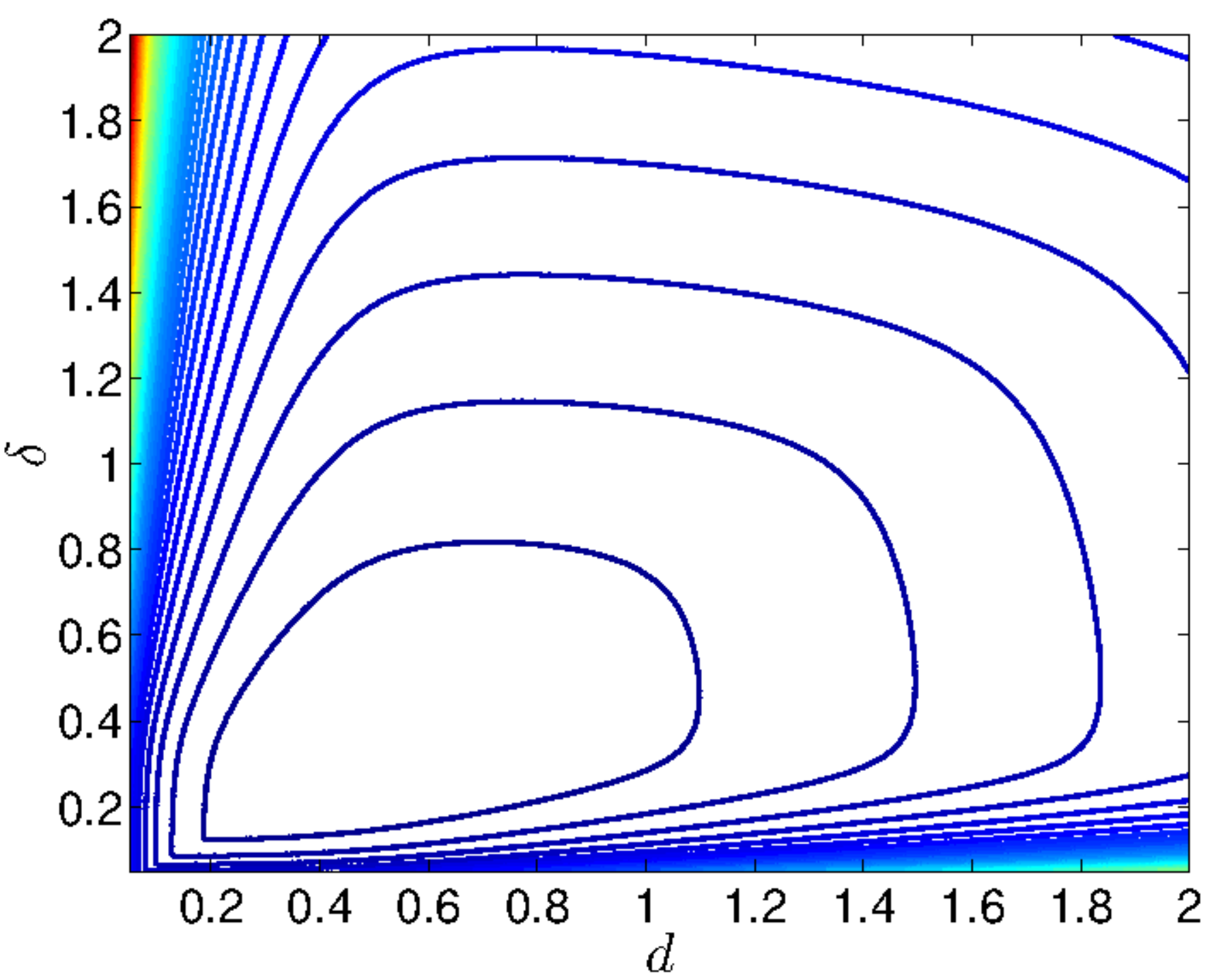,height=6cm, width=8 cm}}

\centerline{\hbox to 9 truecm {\scriptsize (a) \hfill (b)}}

\end{center}
\caption{Objective function \ $\mathcal K(d,\delta)$ \ (a) and the corresponding contour plot (b) of a four-point regular grid design with
\ $\beta=0.2, \ \gamma=0.3$.}
\label{fig4}
\end{figure}
\end{ex}

\subsection{Comparison of equidistant designs}
  \label{subs:subs3.2}
Again, given a fixed design space \ ${\mathcal X}=[a_1,b_1]\times [a_2,b_2], \ a_1<b_1, \ a_2<b_2$, \ we are interested in the effect of refining the directionally equidistant regular grid design by doubling the number of partition intervals in one direction (by symmetry it suffices to deal, for instance, with the first coordinate) or in both coordinate directions.  Further, we also consider the case when doubling of the partition intervals in a given coordinate direction is combined with doubling the corresponding dimension of the design space. Obviously, assumption \ ${\mathcal X}=[0,1]^2$ \ does not violate generality, and we may consider  
designs \ $\boldsymbol \xi_{n,m}$, \ $\boldsymbol \xi_{2n,m}$ \ and \ 
$\boldsymbol \xi_{2n,2m}$, \ on \ $\mathcal X$ \ with \ $\boldsymbol \xi_{n,m}:=\big \{(i/n,j/m): i=0,1, \ldots ,n, \ j=0,1,\ldots ,m\big\}, \ n,m \!\geq\! 2, \ n,m\in{\mathbb N}$, \ design \ $\widetilde{\boldsymbol \xi}_{2n,2m}:=\big \{(i/n,j/m): i=0,1, \ldots ,2n, \ j=0,1,\ldots ,2m\big\}$ \ on \ $\widetilde{\mathcal X}:=[0,2]^2$ \ and  \ $\widehat{\boldsymbol \xi}_{2n,m}:=\big \{(i/n,j/m): i=0,1, \ldots ,2n, \ j=0,1,\ldots ,m\big\}$ \ on \ $\widehat{\mathcal X}:=[0,2]\times [0,1]$. 

\begin{thm}
  \label{2Ddoublim}
For model \eqref{2Dmodel} with covariance structure \eqref{2Doucov}
\begin{equation*}
\lim_{n,m\to\infty}\frac {{\mathcal D}(2n,2m)}{{\mathcal D}(n,m)}=
\lim_{n,m\to\infty}\frac {{\mathcal D}(2n,m)}{{\mathcal D}(n,m)}=
\lim_{n,m\to\infty}\frac {{\mathcal K}(2n,2m)}{{\mathcal K}(n,m)}=
\lim_{n,m\to\infty}\frac {{\mathcal K}(2n,m)}{{\mathcal K}(n,m)}=1,  
\end{equation*}
where \ $\mathcal D(n,m)$ \ and \ $\mathcal K(n,m)$ \ are the values of
the objective functions \eqref{2DDopt} and \eqref{2DKopt}, respectively,
corresponding to the design \ $\boldsymbol \xi_{n,m}$ \ on \ ${\mathcal X}$.
\end{thm} 

Using Theorem \ref{2Ddoublim}, on can formulate a similar conclusion as in the case of OU processes. In particular, after a sufficiently large amount of grid design points there is no need of further refinement of the grid. 

\begin{thm}
  \label{2Ddoubintlim}
For model \eqref{2Dmodel} with covariance structure \eqref{2Doucov}
\begin{equation*}
\lim_{n,m\to\infty}\frac {\widetilde{\mathcal D}(2n,2m)}{{\mathcal D}(n,m)}\!=\!\widetilde D(\beta)\widetilde D(\gamma)
\quad \text{and} \quad
\lim_{n,m\to\infty}\frac {\widehat{\mathcal D}(2n,m)}{{\mathcal D}(n,m)}\!=\!\widetilde D(\beta), \qquad \text{with} \qquad \widetilde D(\beta)\!:=\!\frac{2(\beta\!+\!1)}{\beta\!+\!2}D(\beta),  
\end{equation*}
where \ $D(\beta)$ \ is defined by \eqref{DKform}, whereas \ $\widetilde{\mathcal D}(2n,2m)$ \ and \ $\widehat{\mathcal D}(2n,m)$ \ denote the values of the objective function \eqref{2DDopt} corresponding to designs 
\  $\widetilde{\boldsymbol \xi}_{2n,2m}$ \ and \ $\widehat{\boldsymbol \xi}_{2n,m}$, \ respectively.
\end{thm}

\begin{figure}[t]
\begin{center}
\leavevmode
\hbox{
\epsfig{file=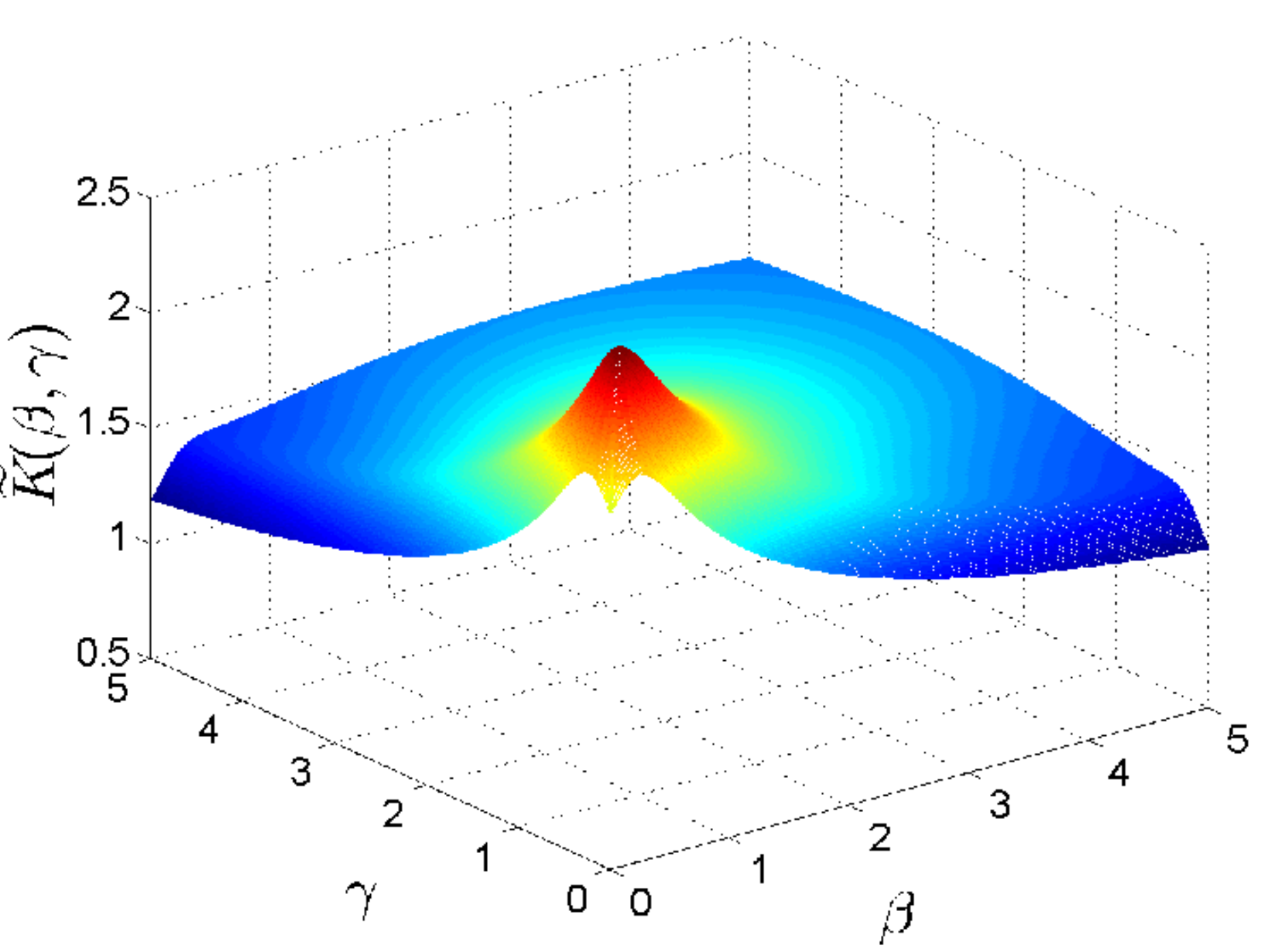,height=6cm, width=8 cm} \quad
\epsfig{file=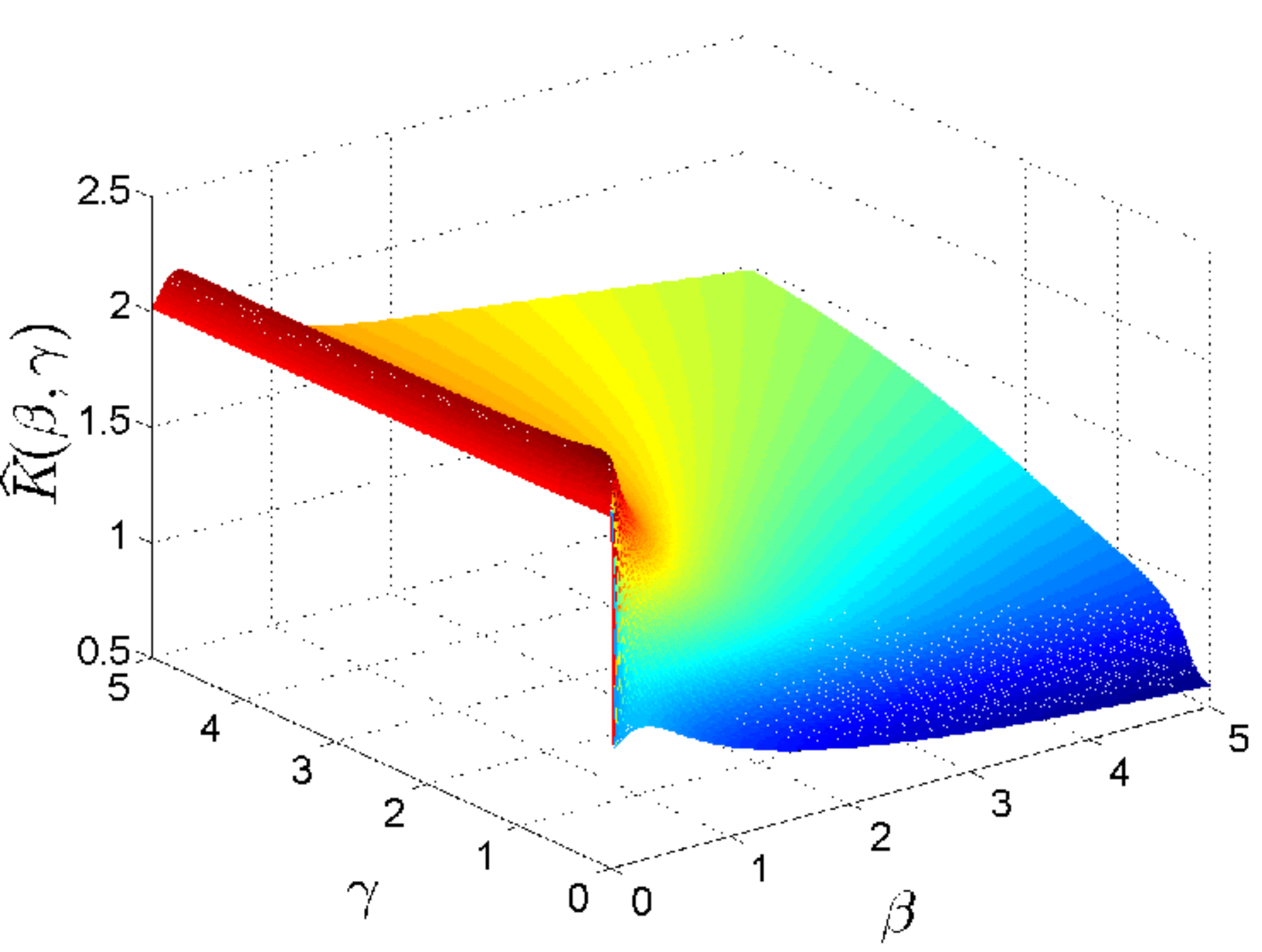,height=6cm, width=8 cm}}

\centerline{\hbox to 9 truecm {\scriptsize (a) \hfill (b)}}

\end{center}
\caption{Limiting functions \ $\widetilde K(\beta,\gamma)$ \ (a) and \ $\widehat K(\beta,\gamma)$ \ (b) of \ $\widetilde{\mathcal K}(2n,2m)/{\mathcal K}(n,m)$  \ and \ $\widehat{\mathcal K}(2n,m)/{\mathcal K}(n,m)$ \ as \ $n,m\to \infty$. }
\label{fig5}
\end{figure}

As \ $\widetilde D(\beta)$ \ is strictly increasing with \ $\lim_{\beta\to 0}\widetilde D(\beta)=2$ \ and \ $\lim_{\beta\to\infty}\widetilde D(\beta)=32$, \ if one has a dense enough directionally equidistant grid of observations, the extension of the design space along a coordinate direction will at least double the information on the unknown regression parameters $(\alpha_0,\alpha_1,\alpha_2)$.

Now, denote by $\widetilde{\mathcal K}(2n,2m)$ \ and \ $\widehat{\mathcal K}(2n,m)$ \  the values of the objective function \eqref{2DKopt} corresponding to designs \  $\widetilde{\boldsymbol \xi}_{2n,2m}$ \ and \ $\widehat{\boldsymbol \xi}_{2n,m}$, \ respectively, and let 
\begin{equation*}
\widetilde K(\beta,\gamma):=\lim_{n,m\to\infty}\frac {\widetilde{\mathcal K}(2n,2m)}{{\mathcal K}(n,m)} \qquad \text{and} \qquad
\widehat K(\beta,\gamma):=\lim_{n,m\to\infty}\frac {\widehat{\mathcal K}(2n,m)}{{\mathcal K}(n,m)}.
\end{equation*}
Due to the very complicated form of the objective function \eqref{2DKopt} one cannot provide feasible expressions for the limiting functions\ $\widetilde K(\beta,\gamma)$ \ and \ $\widehat K(\beta,\gamma)$, \ plotted in Figures \ref{fig5}a and \ref{fig5}b, respectively. In contrast to the D-optimal design, \ $\widehat K(\beta,\gamma)$ \ depends not only on \ $\beta$. \ Moreover, there is a substantial difference compared to the one dimensional case, since \ $\widetilde K(\beta,\gamma)$ \ seems to have a maximum point.

\section{Simulation results}
  \label{sec:sec4}

\begin{figure}[t]
\begin{center}
\leavevmode
\hbox{
\epsfig{file=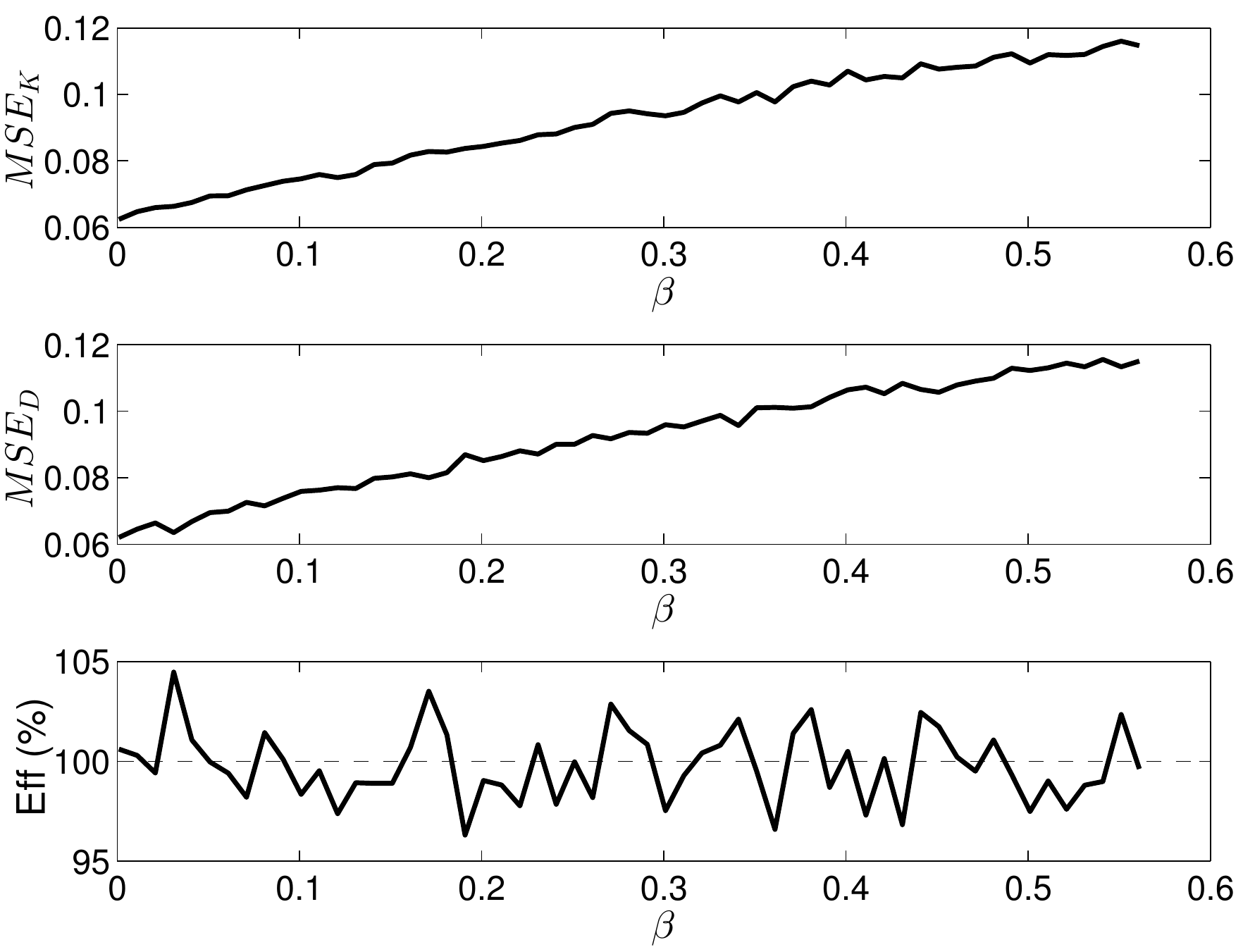,height=6cm, width=8 cm} \quad
\epsfig{file=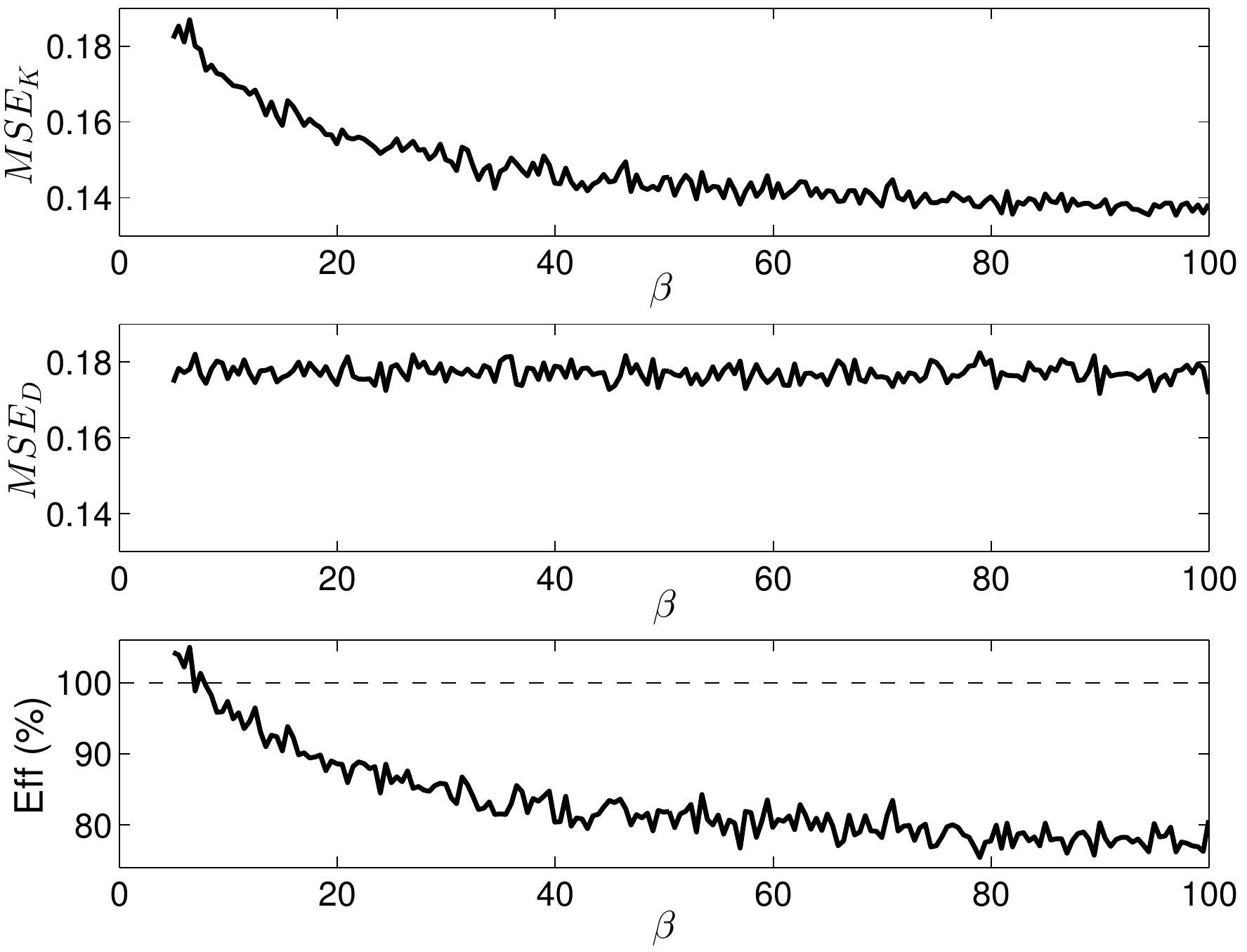,height=6cm, width=8 cm}}

\centerline{\hbox to 9 truecm {\scriptsize (a) \hfill (b)}}

\end{center}
\caption{Average mean squared errors \ $\mse_K$ \ and \ $\mse_D$ \ of GLS estimates of parameters based on three-point restricted K-  and D-optimal designs and relative efficiency \ ($\eff:= {\mse_K}/{\mse_D}\times 100\%$) \ plotted against the parameter \ $\beta$ \ for the  intervals (a) \ $]0,\beta^*[, \ \beta^*\approx 0.5718$, \ and  (b) \  $]\beta^{**},100], \ \beta^{**}\approx 4.9586$, \ respectively.}
\label{fig6}
\end{figure}

To illustrate the differences between K- and D-optimal designs, computer simulations using Matlab R2014a are performed. In general, the driving stationary Ornstein-Uhlenbeck processes and fields of models \eqref{model} and \eqref{2Dmodel}, respectively, can be simulated either 
with the help of discretization \citep[see, e.g.,][]{gil} or using their Karhunen-Lo\`eve expansions based on representations \eqref{ourep} and \eqref{2Dourep} \citep{jb,bs}. However, in our simulation study, due to the small number of locations, it is sufficient to draw samples from the corresponding finite dimensional distributions.

In each of the following examples, \ $10000$ \ independent samples of
the driving Gaussian processes are simulated and the average mean squared errors (MSE) of the generalized least squares (GLS) estimates of the regression parameters based on samples corresponding to different designs are calculated. 

\begin{ex}
  \label{ex4}
Consider the model \eqref{model} with true parameter values \ $\alpha_0=\alpha_1=1$ \ and standard deviation parameter \ $\sigma=1/4$ \ defined on the unit interval \ $[0,1]$, \ and the three-point restricted design \ $\{0,d,1\}$. \ As it has been mentioned in Example \ref{ex1}, the D-optimal design for all \ $\beta>0$ \ is equidistant, whereas K-optimal design exists only for parameter values \ $0<\beta <\beta^*\approx 0.5718$  \ and \ $\beta> \beta^{**}\approx 4.9586$. \ 

Figures \ref{fig6}a and \ref{fig6}b display the average mean squared errors \ $\mse_K$ \ and \ $\mse_D$ \ of the GLS estimates of parameters based on K- and D-optimal designs  plotted against the parameter \ $\beta$ \ together with the relative efficiency 
\begin{equation}
   \label{eff}
\eff:=\frac {\mse_K}{\mse_D}\times 100\%
\end{equation}
for the intervals \ $]0,\beta^*[$ \ and  \ $]\beta^{**},100]$, \ respectively. Observe that for small parameter values the difference in MSE is negligible, whereas for parameters from the upper interval the K-optimal design exhibits a superior overall performance.
\end{ex}

\begin{ex}
  \label{ex5}
Consider the model \eqref{2Dmodel} with true parameter values \ $\alpha_0=\alpha_1=\alpha_1=1$ \ and standard deviation parameter \ $\sigma=1/4$ \ defined on the unit square \ $[0,1]^2$, \ and the nine-point restricted regular grid design \ $\{(0,0),(0,\delta),(0,1),(d,0),(d,\delta),(d,1), (1,0),(1,\delta),(1,1)\}$. \ According to the results of Example \ref{ex2}, for all positive values of \ $\beta$ \ and \ $\gamma$ \ the D-optimal design is directionally equidistant, that is \ $d=\delta=1/2$, \ whereas a non-collapsing K-optimal design exists only in a certain region of the parameter space, see Figure \ref{fig3}. 

\begin{table}[b]
\begin{center}
\hbox{\small
\begin{tabular}{c|rrrrr}
$\beta \setminus \gamma$&0.01&0.03&0.05&0.10&0.15 \\ \hline
0.01&100.88&99.45&97.53&101.60&96.92\\
0.03&99.43&101.54&100.12&99.22&104.38\\
0.05&100.79&100.33&102.00&97.51&97.46\\
0.10&101.37&99.91&99.61&99.86&99.66\\
0.15&100.94&99.51&102.65&98.48&102.48
\end{tabular}
\qquad 
\begin{tabular}{c|rrrrr}
$\beta \setminus \gamma$&10&15&20&25&30 \\ \hline
10&115.33&109.19&105.53&102.89&99.14\\
15&108.18&104.65&96.49&95.25&93.80\\
20&103.93&100.66&95.09&93.61&93.70\\
25&100.50&95.92&94.18&93.20&89.28\\
30&102.91&96.98&94.14&90.04&88.90
\end{tabular}
}
\caption{Relative efficiency $\eff$ (\%) of MSEs of GLS estimates of regression parameters based on nine-point restricted K- and D-optimal regular grid designs for various values of \ $\beta$ \ and \ $\gamma$.}
\label{tab1}
\end{center}
\end{table}
In Table \ref{tab1} the relative efficiencies \eqref{eff} of the MSEs of the GLS estimates of regression parameters based on nine-point restricted K- and D-optimal regular grid designs are reported for different values of \ $\beta$ \ and \ $\gamma$. \ Similar to the temporal case of Example \ref{ex4}, for large parameter values, if a non-collapsing K-optimal design exists, it will outperform the corresponding D-optimal one.
\end{ex}

\section{Conclusions}
  \label{sec:sec5}
We investigate the properties of K-optimal designs for temporal and spatial linear regression models driven by OU processes and sheets, respectively, and highlight the differences compared with the corresponding D-optimal sampling.
We study the problems of existence of K-optimal designs and also investigate the dependence of the two designs on the covariance parameters of the driving processes. This information may be crucial for an experimenter in order to increase efficiency in practical situations. Simulation results display the superiority of restricted K-optimal designs for large covariance parameter values.  

\subsection*{Acknowledgments}
This research was supported by
the J\'anos Bolyai Research Scholarship of the Hungarian Academy of Sciences and by the Hungarian -- Austrian intergovernmental S\&T cooperation program T\'ET\_{}15-1-2016-0046. The author is indebted to Milan Stehl\'\i k for his valuable suggestions and remarks.

\begin{appendix}
\section{Appendix}
  \label{sec:secA}

\subsection{Correlation structure of observations}
   \label{subs:subsA.1}
The correlation matrix of observations \ $\big\{ Y(s_i), \
i=1,2,\ldots ,n\big\}, \ n\geq 2,$ \ of the stochastic process \eqref{model}
equals
\begin{equation*}
C(n)=
     \begin{bmatrix}
        1 &p_1 &p_1p_2  &p_1p_2p_3 &\dots &\dots &\prod_{i=1}^{n-1}p_i \\
        p_1 &1 &p_2 &p_2p_3 &\dots &\dots &\prod_{n=2}^{n-1}p_i \\
        p_1p_2 &p_2 &1 &p_3 &\dots &\dots &\prod_{i=3}^{n-1}p_i \\
        p_1p_2p_3 &p_2p_3 &p_3 &1 &\dots &\dots &\vdots \\
        \vdots &\vdots &\vdots &\vdots &\ddots & &\vdots \\
        \vdots &\vdots &\vdots &\vdots & &\ddots &p_{n-1} \\
        \prod_{i=1}^{n-1}p_i &\prod_{i=2}^{n-1}p_i
        &\prod_{i=3}^{n-1}p_i &\dots &\dots &p_{n-1} &1\\
     \end{bmatrix},
\end{equation*}
where \ $p_i:=\exp(-\beta d_i)$ \  and \ $d_i:=s_{i+1}-s_i, \ i=1,2, \ldots
,n-1$. \ According to the results of \citet{ks}, the inverse of \ 
$C(n)$ \ is given by
\begin{equation*}
  C^{-1}(n)=
     \begin{bmatrix}
        \frac{1}{1-p_1^2} &\frac{p_1}{p_1^2-1} &0 &0 &\dots &\dots &0 \\
        \frac{p_1}{p_1^2-1} &V_{2} &\frac{p_2}{p_2^2-1} &0 &\dots &\dots &0 \\
        0 &\frac{p_2}{p_2^2-1} &V_{3} &\frac{p_3}{p_3^2-1} &\dots &\dots &0 \\
        0 &0 &\frac{p_3}{p_3^2-1} &V_{4} &\dots &\dots &\vdots \\
        \vdots &\vdots &\vdots &\vdots &\ddots & &\vdots \\
        \vdots &\vdots &\vdots &\vdots & &V_{n-1} &\frac{p_{n-1}}{p_{n-1}^2-1} \\
        0 &0 &0 &\dots &\dots &\frac{p_{n-1}}{p_{n-1}^2-1} &\frac{1}{1-p_{n-1}^2}\\
     \end{bmatrix},
\end{equation*}
where \ $V_k:=\frac{1-p_k^2p_{k-1}^2}{(p_k^2-1)(p_{k-1}^2-1)}=\frac
1{1-p_k^2}+ \frac {p_{k-1}^2}{1-p_{k-1}^2}, \ k=2,\dots,n-1$.

\subsection{Proof of Theorem \ref{equiresgen}}
      \label{subs:subsA.2}
Using form \eqref{entry_d} of the entries of the information matrix, a
short calculation shows 
\begin{equation*}
  {\mathcal D}(d)=J_n(\beta d) \frac{n\!-\!1}{\beta^2}F_n(\beta d), 
\end{equation*}
with
\begin{equation}
  \label{detdecomp} 
J_n(d):=\frac {2\!-\!n\!+\!n{\mathrm e}^d}{{\mathrm e}^d\!+\!1}\geq 1
\quad \text{and} \quad 
  F_n(d):=\frac{d^2}{{\mathrm e}^{2d}\!-\!1}\left(
    \frac{n(n\!+\!1)}{12}\big({\mathrm e}^{d}\!-\!1\big)^2 
   \!+\!  \frac{n\!+\!1}2\big({\mathrm e}^{d}\!-\!1\big)\!+\!1\right)\geq 0, 
\end{equation}
where \ $J_n(d)$ \ is a strictly increasing function of
\ $d$. \ Hence, in order to prove the first statement of Theorem
\ref{equiresgen}, it suffices to show that  \ $F_n(d)$ \ is also 
strictly increasing for all integers \ $n\geq 2$. \ This latter property
obviously holds for 
\begin{equation*}
F_2(d)=\frac {d^2}{2(1-{\mathrm e}^{-d})} \qquad \text{and} \qquad
F_3(d)=\frac {d^2}{1-{\mathrm e}^{-2d}},   
\end{equation*}
whereas for \ $n\geq 4$ \ one can consider the decomposition
\begin{equation*}
F_n(d)=F_3(d)G_n(d), \qquad \text{where} \qquad  
 G_n(d):=\frac 1 {{\mathrm
     e}^{2d}}\left(\frac{n(n+1)}{12}\big({\mathrm e}^{d}-1\big)^2 
   +  \frac{n+1}2\big({\mathrm e}^{d}-1\big)+1\right).
\end{equation*}
A short calculation shows that the numerator of \ $G'_n(d)$ \ equals
\begin{equation*}
\frac {{\mathrm e}^{2d}(n-3)}6\big({\mathrm e}^{d}(n+1)-n+2\big )\geq 
\frac {{\mathrm e}^{2d}(n-3)}2>0,
\end{equation*}
which completes the proof of monotonicity of \ ${\mathcal D}(d)$.

Now, \eqref{Ropt} and \eqref{entry_d} imply
\begin{equation}
  \label{Kopt_d}
\mathcal R(d)=\frac {\Big(n({\mathrm e}^{\beta
    d}\!-\!1)^2\big((n\!-\!1)(2n\!-\!1)d^2+6\big)/6 +({\mathrm e}^{\beta
    d}\!-\!1)\big(n(n\!-\!1)d^2+2 \big)+(n\!-\!1)d^2 \Big)^2}
{d^2(n\!-\!1)({\mathrm e}^{\beta 
    d}\!-\!1)\big(n(n\!+\!1)({\mathrm e}^{\beta
    d}\!-\!1)^2/12+(n\!+\!1)({\mathrm e}^{\beta 
    d}\!-\!1)/2+1 \big)\big(n({\mathrm e}^{\beta d}\!-\!1)+2\big)}. 
\end{equation}
After dividing both the numerator and the denominator of the right-hand
side of \eqref{Kopt_d} by \ ${\mathrm e}^{4\beta d}$, \ one can easily
see
\begin{equation}
   \label{Klim_r}
\lim_{d\to\infty} \mathcal R(d)=\lim_{d\to\infty}\frac{n^2\big((n-1)(2n-1)d^2+6\big)^2}{3n^2
  \big(n^2-1\big)d^2}=\infty.    
\end{equation}
In a similar way, taking into account that \ $\lim_{d\to 0}\big({\mathrm e}^d-1\big)/d=1$, \ the division of  both the numerator and the denominator of \ ${\mathcal R}(d)$ \ by \ $d^2({\mathrm e}^{\beta d}-1)$ \ results in
\begin{equation*}
\lim_{d\to 0} \mathcal R(d)=\lim_{d\to 0}\frac {\Big(nd^{1/2}\big((n-1)(2n-1)d^2+6\big)/6 +\big(n(n-1)d^{3/2}+2d^{-1/2} \big)+(n-1)d^{1/2} \Big)^2}{2(n-1)} =\infty,
\end{equation*}
which together with \eqref{Klim_r} implies that \ $\mathcal R(d)$ \ should have at least one global minimum. \proofend

\subsection{Proof of Theorem \ref{eqirestwop}}
      \label{subs:subsA.3}
For \ $n=2$ \ expression \eqref{Kopt_d} simplifies to
\begin{equation*}
\mathcal R(d)=\frac {\big((d^2+2){\mathrm e}^{\beta
    d}-2\big)^2}{d^2\big({\mathrm e}^{2\beta d}-1\big)},
\end{equation*}
and
\begin{equation*}
\mathcal R'(d)=\frac{2\big((d^2\!+\!2){\mathrm e}^{\beta d}\!-\!2\big)}{d^3
\big({\mathrm e}^{2\beta d}\!-\!1\big)^2} R(d) \quad \text{with} \quad
R(d):= (d^2-2){\mathrm e}^{3\beta d}+2(\beta d+1){\mathrm e}^{2\beta
  d} - (\beta d^3+d^2+2\beta d -2){\mathrm e}^{\beta d}-2. 
\end{equation*}
For \ $d\geq 0$ \ the derivative \ $\mathcal R'(d)$ \ 
equals \ $0$ \ if and only if equation \eqref{twopeq} holds, that is \ $R(d)=0$. \
Hence, to complete the proof of Theorem \ref{eqirestwop}, it remains
two show that \eqref{twopeq} has a unique non-negative solution and
this solution is the minimum of \ $\mathcal R(d)$. \  
Now, observe that \ $R(d)$ \ admits the representation
\begin{equation*}
R(d)=d^2R^{(1)}(\beta d)-2R^{(2)}(\beta d),
\end{equation*}
where
\begin{equation*}
R^{(1)}(x):={\mathrm e}^x\big({\mathrm e}^{2x}-x-1\big)>0
\quad \text{and} \quad R^{(2)}(x):=
\Big({\mathrm e}^x\big({\mathrm e}^x-x-1\big)+\big({\mathrm e}^x-1\big)\Big)\big({\mathrm e}^x-1\big)>0, \ \ \text{for} \ \ x>0.
\end{equation*}
In this way
\begin{equation}
  \label{Rroot}
R(d)=0 \qquad \text{if and only if} \qquad \frac {d^2}2=\frac {R^{(2)}(\beta d)}
{R^{(1)}(\beta d)}. 
\end{equation}
A short calculation shows that \ ${R^{(2)}(x)}/ {R^{(1)}(x)}$ \ is strictly monotone increasing and strictly concave. Further, we have 
\begin{equation*}
\lim_{x\to 0} \frac {R^{(2)}(x)} {R^{(1)}(x)}=0 \qquad \text{and} \qquad
\lim_{x\to \infty} \frac {R^{(2)}(x)} {R^{(1)}(x)}=1.
\end{equation*}
Hence, by the convexity of \ $d^2/2$, \ for all \ $\beta >0$ \ the equation on the right hand side of \eqref{Rroot} has a single solution (obviously, depending on \ $\beta$). \ Moreover, as for all \ $\beta >0$, 
\begin{equation*}
\lim_{d\to 0}R(d)=0 \qquad \text{and} \qquad \lim_{d\to \infty}R(d)=\infty,
\end{equation*}
if  \ $R(d)$ \ has a change of sign at this unique root, the change shall be from negative to positive. This means that the solution of \eqref{twopeq} is the unique minimum of \ $\mathcal R(d)$. \ Thus, it remains to show that for all \ $\beta >0$, \ one can find some \ $d>0$ \ such that \ $R(d)<0$. 

Consider the decomposition
\begin{equation*}
R(d)={\mathrm e}^{\beta d} \big(d^2Q^{(1)}(d) -2{\mathrm e}^{\beta
  d}Q^{(2)}(d)\big)-2Q^{(3)}(d), 
\end{equation*}
where
\begin{align*}
&Q^{(1)}(d):={\mathrm e}^{2\beta d}-
\beta^2\Big(\frac{\beta d}3+1\Big) {\mathrm e}^{\beta d} -\beta
d-1, \\[1mm]
Q^{(2)}(d):={\mathrm e}^{\beta d}&\,-
\Big(\frac{\beta^3 d^3}6+\frac{\beta^2d^2}2+\beta d+1\Big), \qquad
Q^{(3)}(d):={\mathrm e}^{\beta d}(\beta d-1)+1.
\end{align*}
Using Taylor series expansions one can easily see that for all \ $x>0$
\begin{equation*}
1-x <{\mathrm e}^{-x} \qquad \text{and} \qquad 1+x+\frac{x^2}2+\frac
{x^3}6 <{\mathrm e}^x,  
\end{equation*}
implying for all \ $d>0$ \ and \ $\beta>0$ \ the positivity of \
$Q^{(3)}(d)$ \ and \ $Q^{(2)}(d)$, \ respectively. Finally,
\begin{equation*}
\lim_{d\to 0} Q^{(1)}(d)=-\beta^2 \qquad \text{and} \qquad \lim_{d\to
  \infty} Q^{(1)}(d)=\infty,
\end{equation*}
thus, for all \ $\beta >0$ \ there exist a \ $d>0$ \ such that \
$Q^{(1)}(d)<0$, \ which completes the proof. \proofend

\subsection{Proof of Theorem \ref{doublim}}
      \label{subs:subsA.4}
Under the settings of the theorem \ ${\mathcal D}(n)= L_1(n+1)L_3(n+1)-L_2^2(n+1)$ \ and \ ${\mathcal K}(n)=g\big(\mathcal R(n)\big)$ \ with \ ${\mathcal R}(n)=\big(L_1(n+1)+L_3(n+1)\big)^2 / \mathcal D(n)$ \
and \ $g(x):=\frac 14 \big(\sqrt{x}+\sqrt{x-4}\big)^2$, \ 
where the expressions for \ $L_1(n+1), \ L_2(n+1)$ \ and \
$L_3(n+1)$ \ can be obtained using \eqref{entry_d} with \
$d=1/n$. \ 
Since for all \ $\beta>0$
\begin{equation*}
 \lim_{n\to\infty} n\big({\mathrm e}^{\beta /n}-1\big)=\beta,
\end{equation*}
one can easily show
\begin{align}
   \label{Llim}
\lim_{n\to\infty}L_1(n+1)=&\lim_{n\to\infty}L_1(2n+1)=\frac{\beta} 2+1,
\qquad \lim_{n\to\infty}L_2(n+1)=\lim_{n\to\infty}L_2(2n+1)=\frac{\beta} 4+\frac 12, \\
&\lim_{n\to\infty}L_3(n+1)=\lim_{n\to\infty}L_3(2n+1)=\frac 1{2\beta}\left(\frac{\beta^2}3+\beta +1\right),  \nonumber
\end{align}
which completes the proof. \proofend

\subsection{Proof of Theorem \ref{doubintlim}}
      \label{subs:subsA.5}

Similar to the proof of Theorem \ref{doublim} we have \ 
$\widetilde{\mathcal D}(2n)= \widetilde  L_1(2n+1)\widetilde L_3(2n+1)-\widetilde L_2^2(2n+1)$ \ and \ $\widetilde{\mathcal K}(2n)=g\big(\widetilde{\mathcal R}(2n)\big)$ \ with \ 
$\widetilde {\mathcal R}(2n)=\big(\widetilde L_1(2n+1)+\widetilde L_3(2n+1)\big)^2 / \widetilde{\mathcal D}(2n)$, \
where now the expressions for \ $\widetilde L_1(2n+1), \ \widetilde L_2(2n+1)$ \ and \ $\widetilde L_3(2n+1)$ \ can be obtained using \eqref{entry_d} for \ $2n+1$ \ design points with \ $d=1/n$. \ Hence, the limits in \eqref{kint_lim} are direct consequences of \eqref{Llim} and 
\begin{equation*}
\lim_{n\to\infty}\widetilde L_1(2n+1)=\beta +1,
\qquad \lim_{n\to\infty}\widetilde L_2(2n+1)=\beta +1 , \qquad
\lim_{n\to\infty}\widetilde L_3(2n+1)=\frac 1{\beta}\left(\frac 43{\beta^2}+2\beta +1\right). 
\end{equation*}
\proofend

\subsection{Proof of Theorem \ref{FIMsheet}}
      \label{subs:subsA.6}

For the regular grid design introduced in Section \ref{sec:sec3}, the covariance matrix \ $C(n,m)$ \ of observations admits the decomposition
\begin{equation*}
C(n,m)=P(n)\otimes Q(m),
\end{equation*}
where \ $P(n)$ \ and \ $Q(m)$ \ are covariance matrices of observations of OU processes with covariance parameters \ $\beta>0$ \ and \ $\gamma>0$ \ in time points \ $s_1<s_2<\ldots < s_n$ \ and \ $t_1<t_2<\ldots < t_m$, respectively (see \citet{bss} or the online supplement of \citet{bsschemo}). In this way, 
\begin{equation}
    \label{Cinvdecomp}
C^{-1}(n,m)=P^{-1}(n)\otimes Q^{-1}(m),
\end{equation}
where for the exact forms of $P^{-1}(n)$ \ and \ $Q^{-1}(m)$ \ see \ref{subs:subsA.1}.
Further,
\begin{equation}
  \label{Gdecomp}
G(n,m)=S(n)\otimes {\mathbf 1}_m^{\top}+{\mathbf 1}_n^{\top}\otimes T(m),
\end{equation}
where \ ${\mathbf 1}_n, \ n\in{\mathbb N},$ \ denotes the column vector of ones of length \ $n$, \ 
\begin{equation*}
S(n):=\begin{bmatrix} 1/2 & 1/2 & \cdots & 1/2 \\ s_1 & s_2 & \cdots &
  s_n \\ 0 & 0 & \cdots & 0 \end{bmatrix} \quad \text{and} \quad
T(m):=\begin{bmatrix} 1/2 & 1/2 & \cdots & 1/2 \\ 0 & 0 & \cdots & 0 \\
t_1 & t_2 & \cdots & t_m \end{bmatrix}.
\end{equation*}
Decompositions \eqref{Cinvdecomp} and \eqref{Gdecomp} and the properties of the Kronecker product imply
\begin{align}
  \label{FIMdecomp}
{\mathcal I}_{\alpha_0,\alpha_1,\alpha_2}(n,m)\!=&\,\big(S(n)P^{-1}(n)S^{\top}(n)\big) \!\otimes\! \big({\mathbf 1}_m^{\top}Q^{-1}(m){\mathbf 1}_m\big)\!+\!
\big(S(n)P^{-1}(n){\mathbf 1}_n\big) \!\otimes\! \big({\mathbf 1}_m^{\top}Q^{-1}(m)T^{\top}(m)\big) \\
&+\!\big({\mathbf 1}_n^{\top}P^{-1}(n)S^{\top}(n)\big) \!\otimes\! \big(T(m)Q^{-1}(m){\mathbf 1}_m\big)\!+\!\big({\mathbf 1}_n^{\top}P^{-1}(n){\mathbf 1}_n\big)\!\otimes\! \big(T(m)Q^{-1}(m)T^{\top}(m)\big). \nonumber
\end{align}
Matrix manipulations, similar to the proof of \eqref{entry}, show
\begin{alignat*}{2}
S(n)P^{-1}(n)S^{\top}(n)&=
\begin{bmatrix} \frac{L_1(n)}4 & \frac{L_2(n)}2 & 0 \\ \frac{L_2(n)}2 & L_3(n) & 0\\ 0 & 0& 0 
\end{bmatrix}, \qquad
T(m)Q^{-1}(m)T^{\top}(m)&&=
\begin{bmatrix} \frac{M_1(m)}4 & 0 & \frac{M_2(m)}2 \\  0 & 0& 0 \\ \frac{M_2(m)}2 & 0 &  M_3(m)  
\end{bmatrix},\\
S(n)P^{-1}(n){\mathbf 1}_n&=\begin{bmatrix} \frac{L_1(n)}2 & L_2(n) & 0 \end{bmatrix}^{\top}, \qquad \quad \,
T(m)Q^{-1}(m){\mathbf 1}_m&&=\begin{bmatrix} \frac{M_1(m)}2 & 0 & M_2(m)\end{bmatrix}^{\top},\\[2mm]
{\mathbf 1}_n^{\top}P^{-1}(n){\mathbf 1}_n&=L_1(m),
\qquad \qquad \qquad \qquad \qquad {\mathbf 1}_m^{\top}Q^{-1}(m){\mathbf 1}_m&&=M_1(m),
\end{alignat*}
which together with \eqref{FIMdecomp} implies \eqref{FIMform}. \proofend

\subsection{Proof of Theorem \ref{2Dequiresgen}}
      \label{subs:subsA.7}

Using \eqref{2DDopt} and representations \eqref{entry_d} and \eqref{entry_m} of the FIM, a short calculation shows
\begin{align}
   \label{Dprod}
{\mathcal D}(d,\delta)&=\left(L_1^2(n)\Big(L_3(n)-\frac {d^2(n-1)^2}4 L_1(n)\Big)\right)\left(M_1^2(m)\Big(M_3(m)-\frac {\delta^2(m-1)^2}4 M_1(m)\Big)\right)\\
&=\left( J_n^2(\beta d) \frac{n-1}{\beta^2}F_n(\beta d)\right)\left( J_m^2(\gamma \delta) \frac{m-1}{\gamma^2}F_m(\gamma \delta)\right), \nonumber
\end{align}
where functions \ $J_n(d)$ \ and \ $F_n(d)$ \ are defined by \eqref{detdecomp}. As both functions are strictly increasing in \ $d$ \ for all integers \ $n\geq 2$, \ decomposition \eqref{Dprod} directly implies the statement of Theorem \ref{2Dequiresgen}. \proofend

\subsection{Proof of Theorem \ref{2Ddoublim}}
      \label{subs:subsA.8}
The statement of the theorem is a direct consequence of \eqref{Llim} and the corresponding limits of \ $M_1(m+1), \ M_2(m+1)$ \ and \ $M_3(m+1)$ \ as \ $m\to \infty$, \ where expressions for \ $M_i(m+1), \ i=1,2,3,$ \ can be obtained using \eqref{entry_m} with \ $\delta=1/m$. \proofend
\end{appendix}


\begin{thebibliography}{100}
\bibitem[Abt and Welch, 1998]{aw} Abt, M. and Welch, W. J. (1998)
  Fisher information and maximum-likelihood estimation of covariance
  parameters in Gaussian stochastic processes. {\em Canad. J. Statist.\/} 
  {\bf 26}, 127--137.

\bibitem[Baldi Antognini and Zagoraiou,  2010]{baz} 
Baldi Antognini, A. and Zagoraiou, M. (2010) Exact optimal designs
for computer experiments via Kriging
metamodelling. {\em J. Statist. Plann. Inference.\/} {\bf 140}, 2607--2617.

\bibitem[Baran {\em et al.\/}, 2003]{bpz} Baran, S., Pap, G. and Zuijlen,
  M.v. (2003) Estimation of the mean
  of stationary and nonstationary Ornstein-Uhlenbeck processes and
  sheets. {\em Comput. Math.  Appl.\/} {\bf 45}, 563--579.

\bibitem[Baran and Sikolya, 2012]{bs}  Baran, S. and Sikolya, K. (2012)
  Parameter estimation in linear regression driven by a Gaussian sheet. 
  {\em Acta Sci. Math. (Szeged)\/} {\bf 78},  689--713.

\bibitem[Baran {\em et al.\/}, 2013]{bss13} Baran, S., Sikolya, K. and 
  Stehl\'\i k, M. (2013) On the optimal designs for prediction of 
  Ornstein-Uhlenbeck sheets. {\em Statist. Probab. Lett.\/} {\bf 83},  
  1580--1587.

\bibitem[Baran {\em et al.\/}, 2014]{bss} Baran, S., Sikolya, K. and
  Stehl\'\i k, M. (2014) Optimal designs for the methane flux in troposphere. 
  arXiv:1404.1839. 

\bibitem[Baran {\em et al.\/}, 2015]{bsschemo} Baran, S., Sikolya, K. and
  Stehl\'\i k, M. (2015) Optimal designs for the methane flux in troposphere. 
  {\em Chemometr. Intell. Lab.\/} {\bf 146}, 407--417. 

\bibitem[Dette {\em et al.\/}, 2015]{dpz} Dette, H., Pepelyshev, A. and 
  Zhigljavsky, A. (2015) Design for linear regression models with 
  correlated errors. In: Dean, A., Morris, M., Stufken, J. and Bingham, D.
  (eds.), {\em Handbook of Design and Analysis of Experiments.\/} Chapman \&
  Hall/CRC, Boca Raton, pp. 237--278.

\bibitem[Dette {\em et al.\/}, 2016]{dpz2} Dette, H., Pepelyshev, A. and 
  Zhigljavsky, A. (2016) Optimal designs in regression with correlated 
  errors. {\em Ann. Statist.\/} {\bf 44}, 113--152.

\bibitem[Gillespie, 1996]{gil} Gillespie, D. T. (1996) Exact numerical 
  simulation of the Ornstein-Uhlenbeck process and its integral. {\em Phys.  
  Rev. E\/} {\bf 54}, 2084--2091.

\bibitem[Hoel, 1958]{hoel58} Hoel, P. G. (1958). Efficiency problems in
  polynomial estimation. {\em Ann. Math. Stat.\/} {\bf 29},
  1134--1145.

\bibitem[Jaimez and Bonnet, 1987]{jb} Jaimez, R. G. and Bonnet,
  M. J. V. (1987) On the Karhunen-Lo\`eve expansion for transformed
  processes. {\em Trabajos Estad\'\i st.\/} {\bf 2},  81--90. 

\bibitem[Kiefer, 1959]{kiefer59} Kiefer, J. (1959) Optimum experimental 
  designs (with discussions). {\em J. R. Statist. Soc. B\/} {\bf 21}, 
    272--319.

\bibitem[Kise\v l\'ak and Stehl\'{\i}k, 2008]{ks} Kise\v l\'ak, J. and
  Stehl\'{\i}k, M. (2008) Equidistant D-optimal designs for parameters
  of Ornstein-Uhlenbeck process. {\em Statist. Probab. Lett.\/} 
  {\bf 78}, 1388--1396. 

\bibitem[Mar\'echal {\em et al.\/}, 2015]{myz} Mar\'echal, P., Ye, J. and 
  Zhou, J. (2015) K-optimal design via semidefinite programming and entropy 
  optimization. {\em Math. Oper. Res.\/} {\bf 40}, 495--512.

\bibitem[M\"uller, 2007]{muller} M\"uller, W. G. (2007) {\em Collecting 
    Spatial Data. Third Edition.\/} Springer, Heidelberg.

\bibitem[M\"uller and Stehl\'\i k, 2004]{ms} M\"uller, W. G. and Stehl\'\i k, 
  M. (2004) An example of D-optimal designs in the case of correlated errors. 
  In: Antoch, J. (ed.), {\em COMPSTAT 2004 -- Proceedings in Computational 
  Statistics.\/} Springer, Heidelberg, pp. 1542--1550.

\bibitem[N\"ather, 1985]{nather} N\"ather, W. (1985) {\em Effective Observation 
    of Random Fields.\/} Teubner Verlagsgesellschaft, Leipzig.

\bibitem[P\'azman, 2007]{pazman} P\'azman, A. (2007)
  Criteria for optimal design of small-sample experiments with
  correlated observations. {\em Kybernetika\/} {\bf 43}, 453--462.

\bibitem[Pukelsheim, 1993]{puk} Pukelsheim, F. (1993) {\em Optimal Design of 
    Experiments.} Wiley, New York.


\bibitem[Rempel and Zhu, 2014]{rz} Rempel, M. F. and Zhou, J. (2014) On exact
 K-optimal designs minimizing the condition number. {\em  Comm. Statist. 
 Theory Methods\/} {\bf 43}, 1114--1131.

\bibitem[Shewry and Wynn, 1987]{swin} 
Shewry, M. C. and Wynn, H. P. (1987) Maximum
  entropy sampling. {\em J. Appl. Stat.\/} {\bf 14},  165--170.

\bibitem[Shorack and Wellner, 1986]{sw} Shorack, G. R. and Wellner,
  J. A. (1986) {\em Empirical Processes with Applications to
    Statistics.\/} Wiley, New York. 

\bibitem[Silvey, 1980]{silvey} Silvey, S. D. (1980) {\em Optimal Design.\/} 
  Chapman \& Hall, London.

\bibitem[Smith, 1961]{smith} Smith, O, K. (1961) Eigenvalues of a symmetric 
  $3\times 3$ matrix. {\em Commun.  ACM\/} {\bf 4}, 168. 

\bibitem[Xia et al., 2006]{xia} Xia, G., Miranda, M. L. and
  Gelfand, A. E. (2006) Approximately optimal spatial design approaches
  for environmental health data. {\em Environmetrics\/} {\bf 17},
  363--385.

\bibitem[Ye and Zhou, 2013]{yz} Ye, J. and Zhou, J. (2013) Minimizing the 
  condition number to construct design points for polynomial regression models.
  {\em Siam. J. Optim.\/} {\bf 23}, 666--686.

\bibitem[Zagoraiou and Baldi Antognini, 2009]{zba09}
  Zagoraiou, M. and Baldi Antognini, A. (2009) Optimal designs for
  parameter estimation of the Ornstein-Uhlenbeck
  process. {\em Appl. Stoch. Models Bus. Ind.\/} {\bf 25}, 583--600.
\end{thebibliography}
\end{document}